\IfFileExists{config.tex}{\def\CfgMode{preprint} 

}{
  \def\CfgMode{review} 

}
\RequirePackage[final]{graphicx}
\documentclass[\CfgMode]{elsarticle}

\usepackage{silence}
\usepackage[utf8]{inputenc}
\usepackage[T1]{fontenc}
\usepackage[english]{babel}
\usepackage[x11names]{xcolor}
\usepackage[colorlinks,allcolors=SteelBlue4,final]{hyperref}
\graphicspath{{./}{figures/}}

\usepackage{amsmath}
\usepackage{amsfonts}
\usepackage{amsthm}
\usepackage{mathdots}
\usepackage{multirow}

\makeatletter

\def\case#1{\csname #1@\csname #1\endcsname\endcsname}

\definecolor{chNewColor}{HTML}{236400}
\definecolor{chOldColor}{HTML}{ff3b3b}

\def\CfgHighlightChanges@true{
  \newcommand{\ch}[1]{\begingroup\color{chNewColor} ##1\endgroup}
  \newcommand{\xch}[2]{%
    {\color{chNewColor} ##1}%
    {\color{chOldColor}\scriptsize (##2)}%
  }
}
\def\CfgHighlightChanges@false{
  \newcommand{\ch}[1]{##1}
  \newcommand{\xch}[2]{##1}
}
\case{CfgHighlightChanges}

\def\CfgShowComments@true{
  \usepackage{todobox}

  \usepackage[textwidth=\textwidth,right=22em,marginparwidth=20em]{geometry}
}
\def\CfgShowComments@false{
  \newcommand{\todo}[2][]{}
}
\case{CfgShowComments}

\def\CfgLineNumbers@true{
  \usepackage{vruler}
  \AtBeginDocument{\expandafter\setvruler\expandafter[\the\baselineskip]}
}
\def\CfgLineNumbers@false{

}
\case{CfgLineNumbers}

\makeatother

\usepackage{newtxtext}
\usepackage[bigdelims]{newtxmath}
\usepackage{booktabs}
\usepackage{algorithm}
\usepackage{algpseudocode}
\usepackage{acro}
\acsetup{first-style=footnote}
\usepackage[per-mode=fraction]{siunitx}
\usepackage{subfigure}

\theoremstyle{plain}
\newtheorem{theorem}{Theorem}

\newtheorem{proposition}[theorem]{Proposition}

\theoremstyle{definition}

\theoremstyle{remark}

\DeclareAcronym{ODE}{short=ODE,long={ordinary differential equation}}
\DeclareAcronym{PDE}{short=PDE,long=partial differential equation}
\DeclareAcronym{IVP}{short=IVP,long=initial value problem}
\DeclareAcronym{FFT}{short=FFT,long=fast Fourier transform}
\DeclareAcronym{CF}{short=CF,long=Carathéodory-Fejér}
\DeclareAcronym{REXI}{
  short=REXI,
  long=rational approximation of exponential integrators}
\DeclareAcronym{CPU}{short=CPU,long=Central Processing Unit}
\DeclareAcronym{MPI}{short=MPI,long=Message Passing Interface}
\DeclareAcronym{GMRES}{short=GMRES, long=Generalized Minimal Residual}
\DeclareAcronym{JURECA}{
  short=JURECA,long=Jülich Research on Exascale Cluster Architectures}

\newcommand{\iu}{\mathrm{i}}
\newcommand{\df}[1]{\mathop{\mathrm{d}#1}}

\newcommand{\dotp}[2]{\langle #1, #2 \rangle}

\def\defnumbers#1#2{
  \expandafter\gdef\csname #1N\endcsname {#2}
  \expandafter\gdef\csname #1V\endcsname ##1{ {#2}^{##1}}
  \expandafter\gdef\csname #1M\endcsname ##1##2{ {#2}^{##1 \times ##2} }
}
\defnumbers{R}{\mathbb{R}}
\defnumbers{C}{\mathbb{C}}
\defnumbers{N}{\mathbb{N}}

\DeclareMathOperator{\diag}{diag}

\DeclareMathOperator{\cond}{cond}
\DeclareMathOperator{\sr}{sr}

\renewcommand{\Re}{\mathop{\mathrm{Re}}}
\renewcommand{\Im}{\mathop{\mathrm{Im}}}
\renewcommand{\vec}[1]{\mathbf{#1}}

\newcommand{\PZ}{\phantom{0}}
\newcommand{\PZZ}{\phantom{00}}

\begin{document}

\begin{frontmatter}

\title{Time-parallel simulation of the Schrödinger Equation}

\author[fzj]{Hannah Rittich\corref{mycorrespondingauthor}}
\cortext[mycorrespondingauthor]{Corresponding author}
\ead{h.rittich@fz-juelich.de}

\author[fzj]{Robert Speck}
\ead{r.speck@fz-juelich.de}

\address[fzj]{Jülich Supercomputing Centre, Forschungszentrum Jülich GmbH, 52425 Jülich, Germany}

\begin{abstract}
The numerical simulation of the time-dependent Schrödinger equation for quantum systems is a very active research topic.
Yet, resolving the solution sufficiently in space and time is challenging and mandates the use of modern high-performance computing systems.
While classical parallelization techniques in space can reduce the runtime per time step, novel parallel-in-time integrators expose parallelism in the temporal domain.
They work, however, not very well for wave-type problems such as the Schrödinger equation.
One notable exception is the rational approximation of exponential integrators.
In this paper we derive an efficient variant of this approach suitable for the complex-valued Schrödinger equation.
Using the Faber-Carathéodory-Fejér approximation, this variant is already a fast serial and in particular an efficient time-parallel integrator.
It can be used to augment classical parallelization in space and we show the efficiency and effectiveness of our method along the lines of two challenging, realistic examples.
\end{abstract}

\begin{keyword}
Schrödinger equation \sep parallel-in-time \sep rational approximation of exponential integrators \sep parallel across the method \sep Faber-Carathéodory-Fejér approximation
\end{keyword}

\end{frontmatter}


\section{Introduction}

The time-dependent, non-relativistic Schrödinger equation \citep{Schroedinger1926} is a complex-valued linear partial differential equation (\acs{PDE}) that describes the time-evolution of a quantum system.
Being able to predict the behavior of a quantum system is important for many applications.
Without an analytical, tractable solution, however, numerical methods are needed to evaluate the solution of the \ac{PDE}.
Interest in simulating the time-dependent Schrödinger equation started in
the end of the 1950s \citep{MazurRubin1959}. With the availability of
sufficiently powerful computers, these simulations became increasingly popular
for the investigation of molecular structures around 1970
\citep{McCulloughWyatt1969,McCulloughWyatt1971a,ParkEtAl1970}
and are still relevant today \citep{BalintKurti2010,Li2019,Ullrich2016}.

For this work we restrict ourselves to the single-particle
Schrödinger equation in $d$ dimensions,
\begin{equation}
  \label{eq:schroedinger}
  \iu \hbar \frac{\partial}{\partial t} \psi(\vec{r}, t)
  = \left(-\frac{\hbar^2}{2m} \Delta + V(\vec{r}) \right)
    \psi(\mathbf{r}, t)
  \,,
\end{equation}
where $\hbar$ is the reduced Planck constant, $m$ the mass of the particle,
$\Delta$ the Laplace operator,
$V: \RV{d} \to \RN$ the potential, and
$\psi: \mathbb{R}^d \times \mathbb{R} \to \mathbb{C}$ is the unknown
wave function $\psi$.
Given an initial wave function $\psi_0$ at some time $t$, the Schrödinger
equation can be used to compute the wave function at any later time.
The function $\psi$ encodes the probability distribution of the
position and momentum of the particle.
More precisely, the probability density of the position of the particle at
time $t$ is given by $| \psi(\vec{r},t) |^2$, while
the momentum of the particle is, loosely speaking, encoded in the wave-length
of $\psi$ via the de Broglie relation $\lambda = h/p$.
For more details see, e.g., \citep{Messiah1999,Shankar2014,Griffiths2017}.

The Schrödinger equation is defined on an unbounded domain, which causes
problems for numerical computations. Hence, we introduce a finite, but sufficiently large, domain
$\Omega \subseteq \RV{d}$. We then demand that $\psi$ fulfills the Schrödinger
equation for $\vec{r} \in \Omega$ and conforms to zero Dirichlet boundary
conditions, i.e., we require that
\begin{equation*}
  \psi(\vec{r}, t) = 0
  \quad\text{for}\quad
  \vec{r} \in \partial\Omega
  \,.
\end{equation*}
These boundary conditions imply that the particle leaves the domain $\Omega$
with zero probability. If $\Omega$ is large enough, this is a reasonable
assumption and does not change the outcome of the simulation.

In order to perform such a simulation, the continuous Schrödinger equation has to be discretized both in space and time.
Depending on the dimension, the smoothness of the solution and the dynamics of the system, the resulting numerical method may require fine and advanced discretization schemes to resolve the solution adequately and over long time-scales.
This mandates the application of parallel numerical algorithms on high-performance computing systems.

Classical parallelization techniques primarily target the spatial domain and are very successful in reducing the time-to-solution per time step.
However, this approach can neither mitigate the need for a better resolution in time nor can it scale indefinitely for a fixed-size problem.
One promising remedy is the application of parallel-in-time integration techniques.
They expose parallelism also in the temporal domain, either within each time step,
referred to as \emph{parallelization across the method}, or by computing multiple time steps simultaneously, referred to as \emph{parallelization across the steps}~\citep{Burrage1997}.

Parallel-in-time methods have been applied to a multitude of problems, ranging from reaction-diffusion systems \citep{schoebel2019pfasster} and a kinematic dynamo \citep{ClarkeEtAl_ParallelInTimeIntegrationOfKinematicDynamos} to eddy current problems \citep{FriedhoffEtAl_MultigridReductionInTimeForEddyCurrentProblems}, fusion plasma simulations \citep{SamaddarEtAl2019} as well as power systems \citep{SchroderEtAl2017} and robotics \citep{agboh2019combining}, to name just a few very recent ones.
For further reading, we refer to the comprehensive list of references that
is provided on the community website on parallel-in-time integration\footnote{\url{https://parallel-in-time.org}}.

Yet, many of these approaches fail for wave-type problems, which includes the
non-relativistic Schrödinger equation, we are interested in.
For this class of problems, only very specialized and often purely theoretical ideas exist.
One promising one, which has indeed shown its practical relevance, is the \emph{rational approximation of exponential integrators} (\acs{REXI}).
While it has been well known that certain forms of rational approximations can
be used to compute the matrix exponential in parallel
\citep{TrefethenWeideman_TheExponentiallyConvergentTrapezoidalRule,
       HaleEtAl_ComputingPowerLogAndRelatedMatrixFunctionsByContourIntegrals},
it has been first applied for the construction of a parallel-in-time
solver for wave-type problems in \citep{HautEtAl2016}.
This method targets linear problems and forms a parallelizable approximation of the exponential matrix function using rational functions, which can then be used to approximate the solution of the linear \ac{PDE}.
The approximation is designed in a way such that its evaluation consists mainly of the computation of a sum, where the computation of each summand is expensive.
The benefit of this structure is that each individual summand in the approximation can be computed in parallel.
It can thus be classified as parallel across the method, although its approach allows to take much larger time steps as more classical methods like Crank-Nicolson.
\ac{REXI} has been successfully applied to shallow-water equation on the
rotating sphere
\citep{SchreiberEtAl_ExponentialIntegratorsWithParallelInTimeRationalApproximationsForTheShallowWaterEquationsOnTheRotatingSphere}
and to linear oscillatory problems \citep{HautEtAl2016,SchreiberEtAl2018}, making parallel-in-time integration possible even for these challenging problems.

The rational approximation chosen in the original \ac{REXI} approach presented in \citep{HautEtAl2016}, however, involves taking the real part of a complex quantity.
While the method can still be applied for complex-valued problems such as the Schrödinger equation, it becomes significantly more expensive.

In this paper, we therefore present a variant of the \ac{REXI} method specifically targeted toward complex-valued problems.
We use a variation of the Faber-Carathéodory-Fejér (Faber-\acs{CF}) approximation together with a conformal Riemann mapping, which is tailored for the purely imaginary eigenvalues of the semi-discretized Schrödinger equation.
This approach reduces the cost of the \ac{REXI} method substantially, i.e., fewer summands are needed for the rational approximation, thereby increasing the ratio of accuracy per parallel task.
For a given accuracy, this method imposes a restriction on the time-step size, which is also discussed in this paper.
We note that this restriction is inherent to the \ac{REXI} approach itself and needs to be considered for the original version as well.

We begin by briefly explaining the finite element discretization in space (Section~\ref{sec:fem}) which leads to a system of ordinary differential equations (\acsp{ODE}) that needs to be solved.
Then, we discuss how to solve this system by approximating a certain matrix exponential (Section~\ref{sec:time}) and how this computation can be performed in parallel (Section~\ref{sec:rexi}).
For this approximation we need to find a suitable rational approximation to the exponential function, which we construct by using the Faber-Carathéodory-Fejér method (Section~\ref{sec:faber-cf-full}).
Finally, we apply the method to two challenging, realistic problems, namely the quantum tunneling and the double-slit experiment, analyzing the performance of the method (Section~\ref{sec:num}) and finally discuss the applicability of the numerical method (Section~\ref{sec:outro}) beyond the Schrödinger equation.


\section{Space Discretization}%
\label{sec:fem}%
To simulate the Schrödinger equation, we need an
appropriate discretization of the equation. We start by applying the
method of lines approach to turn the \ac{PDE} into a system of \acp{ODE},
by applying the finite element method
\citep{JohNumerical2009,BrennerScott2000,Braess2007}
to the spatial part of the \ac{PDE}.

The finite element discretization is based on the weak formulation of the
\ac{PDE}.
For a domain $\Omega$, let $H^1(\Omega) \subset \CV{\Omega}$ be the Sobolev
space of order one and let $\mathring{H}^1(\Omega)$ be the subset of $H^1(\Omega)$
that consists of functions whose trace vanishes, i.e., $\mathring{H}^1(\Omega)
:= \left\{ u \in H^1(\Omega) : u |_{\partial\Omega} = 0 \right\}$,
where $u|_{\partial\Omega}$ is the trace of $u$
\citep{BrennerScott2000}.
The weak formulation of the Schrödinger equation \eqref{eq:schroedinger}
is to find $\psi$ such that
\begin{equation*}
  \text{$\psi(t, \cdot) \in \mathring{H}^1(\Omega)$}
  \quad\text{and}\quad
  \iu b\big(\phi, \tfrac{\partial\psi}{\partial t}(t)\big) = a(\phi, \psi(t))
  \quad\text{for all}\quad
  \phi \in \mathring{H}^1(\Omega)
  \,,
\end{equation*}
where $a$ and $b$ are the bilinear forms given by
\begin{equation*}
  a(\phi, \psi) :=
    \int \left(\tfrac{\hbar^2}{2m} \nabla \phi \cdot \nabla \psi
               + V \phi \psi \right) \df{\mathbf{r}}
  \quad\text{and}\quad
  b\big(\phi, \tfrac{\partial\psi}{\partial t}\big) := \hbar
     \int \phi \tfrac{\partial \psi}{\partial t} \df{\mathbf{r}}
\end{equation*}
\citep{AraiEtAl1976,KanesakaEtAl1978}.

We can turn the weak formulation into a discrete problem using the famous
Ritz approach \citep{Ritz1909}, which approximates the solution of the weak
form of the \ac{PDE}. We select a suitable subspace
$V_h \subseteq \mathring{H}^1$ and replace $\mathring{H}^1$ by $V_h$ in the
weak formulation. In our case, we choose
$V_h = \{ u \in \mathring{H}^1(\Omega) : u |_T \in \mathcal{P}_2
  \text{ for all } T \in \mathcal{T} \}$,
where $\mathcal{T}$ is a triangulation of the domain $\Omega$.
More precisely, we use Lagrange finite elements of order $2$ on each
triangle \citep{Braess2007,BrennerScott2000}.
Introducing basis vectors $\chi_1, \dots, \chi_n$ for $V_h$,
we can write the modified weak formulation as
\begin{equation}
  \label{eq:ode-system}
  \iu B \tfrac{\partial u}{\partial t}(t)
  = A u(t) \,,
\end{equation}
where $A, B \in \RM{n}{n}$ with $A_{jk} = a(\chi_j,\chi_k)$, and
$B_{jk} = b(\chi_j, \chi_k)$. The vector $u \in \CV{n}$ contains the basis
coefficients of the approximation of the solution.
Together with suitable initial conditions $u(0) = u_0$, this system of
\acp{ODE} defines the initial value problem (\acs{IVP}) that we intend to solve
using an efficient, parallel-in-time integrator.


\section{Time Discretization}
\label{sec:time}

Some of the most efficient time integration methods for the Schrödinger
equation are based on the approximation of the matrix exponential
\citep{LeforestierEtAl1991}. Classical time integration schemes require the
step size of the method to be a fraction of the shortest wave-length that
is present in the problem. Methods based on the computation of the matrix
exponential usually do not have this restriction and can in principle use
much larger step sizes.

The matrix exponential can be used to compute the solution of the \ac{IVP}.
Since the matrices $A$ and $B$ do not depend on the time variable $t$,
the solution $u(\tau)$ for $\tau > 0$ of the \ac{IVP}~\eqref{eq:ode-system}
is given by
\begin{equation}
  \label{eq:matrix-exp-solution}
  u(\tau) = \exp(\tau M) \, u_0
  \quad \text{with} \quad
  M := (\iu B)^{-1} A
  \,,
\end{equation}
where $\exp$ is the matrix exponential \citep{Bellman1997,LiesenMehrmann2015}.
Thus, one way to solve the \ac{IVP} is to compute an approximation to the
product of the vector $u_0$ and the matrix exponential of $M$.

\subsection{Rational Approximation of Exponential Integrators} 
\label{sec:rexi}

We want to use rational approximations to compute the matrix exponential
numerically. There are various ways to compute the exponential of a matrix
\citep{MolerVanLoan1978,MolerVanLoan2003}, however,
we are interested in methods that use rational approximations, because
these methods can be constructed in a way that allows for parallelizing the
time integration scheme itself, increasing the parallelism of the overall
solution process
\citep{HautEtAl2016,SchreiberEtAl2018,
SchreiberEtAl_ExponentialIntegratorsWithParallelInTimeRationalApproximationsForTheShallowWaterEquationsOnTheRotatingSphere}.

It can be shown that the matrix $\tau M = \tau \, (\iu B)^{-1} A$ is
diagonalizable with purely imaginary eigenvalues.
Hence, for simplicity, we restrict ourselves in the following to the
computation of exponentials of matrices that have these properties.
The method, however, can be applied in cases where these two assumptions
do not hold.

The matrix exponential is a special case of a matrix function, which is a way
to extend a scalar function $f: \CN \to \CN$ to the set of
matrices, i.e., to a function $f: \CM{n}{n} \to \CM{n}{n}$
\citep{Higham2008}.
If $G$ is diagonalizable, i.e., $G = X \Omega X^{-1}$,
$\Omega = \diag(\omega_1, \dots, \omega_n)$,
the matrix function of $f$ is given by
$f(G) = X \diag(f(\omega_1), \dots, f(\omega_n))\,X^{-1}$.
Diagonalizing a large matrix, however, is computationally expensive, and thus
this formula is usually not useful for computing the function of a matrix.

We can reduce the computational costs by replacing the direct computation
of the matrix function by a suitable approximation.
By using the diagonalization of $G$, we see that
if $\tilde{f}$ is a function which approximates $f$ in the eigenvalues of
$G$ then the matrix function $\tilde{f}(G)$ is close to $f(G)$.
Thus, if $\tilde{f}(G)$ is cheap to compute, we have a practical way
for evaluating the matrix function $f(G)$ numerically.

Matrix functions of rational functions can be computed without the need of
explicitly computing the diagonalization of the matrix.
Assume that $r = p/q$, and $p,q$ are polynomials
such that $r$ approximates $f$ in the eigenvalues of $G$.
Computing $r(G)$ to approximate $f(G)$ is a feasible approach for the
numerical evaluation of the matrix exponential on its own,
however, by making an additional assumption, we can derive a
time-stepping scheme that intrinsically allows for the simultaneous execution
of certain parts of the computation.

Assume that $\deg{p} < \deg{q}$ and that the roots of $p$ are
distinct. In this case, we can use the partial fraction decomposition to
obtain that
\begin{equation}
  \label{eq:partial-fraction-decomposition}
  r(z) = \sum_{j = 1}^K \frac{\beta_j}{z - \sigma_j}
\end{equation}
for proper shifts $\sigma_j\in\CN$ and coefficients $\beta_j\in\CN$,
$j = 1, \dots, K$, $K \in \NN$, and the corresponding matrix function
\begin{equation*}
  r(G) = \sum_{j = 1}^K \beta_j (G - \sigma_j I)^{-1}
  \,,
\end{equation*}
which approximates $f(G)$.

Using the rational approximation of the matrix function $f(G)$, we can
define the \ac{REXI} time stepping scheme.
Let $f = \exp$, $G = \tau M$, and $r$ a rational approximation as discussed
above, i.e., $r(\tau M) \approx \exp(\tau M)$.
Then, the exponential formula for the solution of the \ac{ODE}
\eqref{eq:matrix-exp-solution} implies that $u(\tau) \approx r(\tau M) u_0$.
Hence, we define one time step of the \ac{REXI} method by
$u_1 = r(\tau M)\, u_0$, which can be computed by
\begin{equation}
  \label{eq:rational-evaluation}
  \begin{aligned}
    u_1
    &= \sum_{j = 1}^K \beta_j (\tau (\iu B)^{-1} A - \sigma_j I)^{-1} u_0
    = \sum_{j = 1}^K \beta_j ((\iu B)^{-1} (\tau A - \sigma_j \, \iu B))^{-1} u_0 \\
    &= \sum_{j = 1}^K \beta_j (\tau A - \sigma_j \, \iu B))^{-1}\, (\iu B)\, u_0
  \,,
  \end{aligned}
\end{equation}
where we used the definition of $M$.

The benefit of computing matrix exponential times vector by evaluating the
rational approximation via \eqref{eq:rational-evaluation}
is that the evaluation of this approximation can be readily parallelized,
because each summand can be evaluated independently.
Thus, each of the $K$ linear systems can be solved independently, using $K$ different parallel tasks.
We refer to this particular splitting of the
computation into tasks as \emph{time-parallelization},
because it uses only properties that are inherent to the time-stepping
scheme itself.
Note that each of these $K$ temporal tasks can be parallelized themselves,
since they involve a set of vector and matrix routines, which can be executed
in parallel as well.
We call this second splitting \emph{space-parallelization}, because the
vectors and matrices describe the spatial dimension of the problem.
Using time-parallelization into $P_\mathrm{t}$ tasks and then applying
space-parallelization into $P_\mathrm{s}$ sub-tasks to each of the
temporal tasks yields $P_\mathrm{t} \cdot P_\mathrm{s}$ sub-tasks that
can be executed simultaneously.

To be able to implement and apply this method, it remains to derive proper shifts $\sigma_j$
and coefficients $\beta_j$ (for $j = 1, \dots, K$).
Note that, in general, we only need to compute these shifts and coefficients
once, because they do not depend on the initial values or the time step.
In the following, we will describe the derivation of these parameters in detail using the Faber-Cara\-théo\-dory-Fejér (Faber-CF) approximation, introduced in \citep{EllFaber1983}.
The intention here is to allow interested readers to comprehend and reproduce
the steps necessary to obtain the $\sigma_j$ and $\beta_j$ and thus the full
algorithm.

We point out that the use of the Faber-CF approximation is a key difference to the \ac{REXI}
approach in \citep{HautEtAl2016} and \citep{SchreiberEtAl2018}.
There, an approximation of the form
\begin{equation*}
  e^{\iu x}
  \approx \Re \left( \sum_{j=1}^K \frac{\gamma_j}{\iu x + \mu_j}
  \right)
\end{equation*}
for certain $\gamma_j, \mu_j \in \CN$ is used.
If all eigenvalues of $M$ are purely imaginary, as in our case, and all eigenvectors can be
chosen to be real, then
\begin{equation*}
  e^{\tau M}
  \approx \Re \left( \sum_{j=1}^K \gamma_j (\tau M + \mu_j I)^{-1} \right)
  \,,
\end{equation*}
where $\Re$ denotes the element-wise real-part of the matrix.

The problem when applying this approximation to the Schrödinger equation is
that we want to compute $e^{\tau M}\, u_0$ where $u_0$ has complex entries
without explicitly computing the matrix $e^{\tau M}$.
If $u_0$ would be real, then $u_0$ could just be moved inside the computation
of the real part. Since $u_0$ is complex, however, we need to compute
\begin{equation*}
  e^{\tau M}\, u_0
  \approx
    \Re \left( \sum_{j=1}^K \gamma_j (\tau M + \mu_j I)^{-1} \Re(u_0) \right)
    + \iu \Re \left( \sum_{j=1}^K \gamma_j (\tau M + \mu_j I)^{-1} \Im(u_0) \right)
  \,,
\end{equation*}
which is twice as much work as in the real case.
The Faber-\ac{CF} approximation that we use does not have this drawback.

Furthermore, the shifts used in the method derived in
\citep{HautEtAl2016,SchreiberEtAl2018} form conjugate pairs. Using properties
of the real numbers, the method only needs to solve one linear system for each
conjugate pair. Since the matrix $M$ of the discretization of the
Schrödinger equation has complex entries and the right-hand side of the
linear systems are complex valued, such an simplification is not possible in
the setting we consider in this paper.

There exists another difference between the two approaches.
The Faber-\ac{CF} approximation computes the approximation essentially in
one step, while the method in \citep{HautEtAl2016} involves a two step
approximation. First, a rational approximation to a Gaussian function is
constructed. Then, this approximation is used, to approximate the function
$e^{\iu x}$.
This procedure has the benefit that it is easy to compute approximations
that are accurate over large intervals and thus allow for the large time steps
(see Section~\ref{sec:step-sizes}). In our experience, however, using the
same accuracy and same approximation interval, the Faber-\ac{CF} approximation
requires fewer poles and therefore fewer linear systems to solve, as
detailed in Section~\ref{sec:step-sizes} below.

\subsection{Faber-Carathéodory-Fejér Approximation} 
\label{sec:faber-cf-full}

The Faber-Carathéodory-Fejér approximation is based on the Carathéodory-Fejér (\acs{CF})
approximation introduced in \citep{Trefethen1981b}.
The latter computes an approximation to holomorphic functions on the
unit disc. The former uses the Faber transform, to generalize the
Carathéodory-Fejér approximation to almost arbitrary approximation domains.

\subsubsection{The Carathéodory-Fejér Approximation}

The \ac{CF} approximation
is a rational approximation to a holomorphic function on the unit disc
$D := \{ z \in \CN : | z | \le 1 \}$.
The resulting rational approximations are only close to the best
approximation, but easier to obtain.

Let us start by introducing the following notation.
Let $R_{mn}$ be the set of rational functions $r(z) = \frac{p(z)}{q(z)}$
with $\deg(p) \le m$ and $\deg(q) \le n$ that are holomorphic in
$D$. Furthermore, let $S := \{ z : | z | = 1\}$ denote the unit circle.
We define the uniform norm of a complex valued function $u$ on the unit disc
by $\| u \|_D := \sup \{ | u(z) | : z \in D \}$ and
the uniform norm of a complex valued function $u$ on the unit circle by
$\| u \|_S := \sup \{ | u(z) | : z \in S\}$.

With this notation at hand, our next step is to simplify the problem.
Assume that $r \in R_{mn}$ is an approximation to a
function $f$ which is holomorphic on $D$. In this case, the error
$e(z) := f(z) - r(z)$
is also holomorphic on $D$.
We want that the size of the error to be as small as possible, i.e.,
we want $\| e \|_D$ to be small. Since $e$ is holomorphic on $D$, its maximum
is located on the boundary of $D$.
Thus, to minimize $\| e \|_D$ we just have to minimize $\| e \|_S$,
i.e., we just have to find a rational function that approximates $f$ well
on the sphere $S$.

It is difficult to find the best approximation to $f$ in $R_{mn}$.
The key idea of the \ac{CF} approximation is to find the best approximation
to $f$ in a larger space $\widetilde{R}_{mn}$ with respect to the
$\| \cdot \|_S$ norm and then approximate the best approximation from
$\widetilde{R}_{mn}$ with a function from $R_{mn}$.
The space $\widetilde{R}_{mn}$ is defined as follows.
Let $G$ be the set of functions that are analytic and bounded in
$1 < | z | \le \infty$ and zero at $z = \infty$. Then define
$\widetilde{R}_{nn} := R_{nn} + G$ and $\widetilde{R}_{mn} := z^{m-n} \widetilde{R}_{nn}$.
One can show that
the space $\widetilde{R}_{mn}$ consists of the functions
\begin{equation}
  \label{eq:r-tilde-form}
  \tilde{r}(z) = \frac{\sum_{k = -\infty}^m d_k z^k}
                      {\sum_{k = 0}^n e_k z^k}
  \,,
\end{equation}
where the poles of the numerator lie inside the unit disc and
the roots of the denominator lie outside the unit disc.

Once we have obtained the best approximation $\tilde{r}^*\in\widetilde{R}_{mn}$ in the form of
\eqref{eq:r-tilde-form} we can use it to find an approximation in
$R_{mn}$ that is close to $\tilde{r}^*$---the \ac{CF} approximation.
Consider the asymptotic analysis of approximating a function
$z \mapsto f(\epsilon z)$ for $\epsilon \to 0$, $\epsilon > 0$ where
$f$ is smooth.
In \citep{Trefethen1981b} it was shown that for small enough $\epsilon$, the best
approximation $\tilde{r}^*$ gets arbitrarily close to a rational function.
This behavior motivates the construction of the \ac{CF} approximation
$r^\mathrm{cf}$:
we compute $\tilde{r}^*$ in the form of \eqref{eq:r-tilde-form}
and discard the summands with negative indices from the numerator, i.e.,
\begin{equation*}
  r^\mathrm{cf}(z) = \frac{\sum_{k = 0}^m d_k z^k}
                          {\sum_{k = 0}^n e_k z^k}
  \,.
\end{equation*}

We thus need to find the best approximation $\tilde{r}^*$ of $f$ in
$\widetilde{R}_{mn}$.
First of all note that $f$ can be written in Maclaurin series form,
because $f$ is holomorphic. In case the Maclaurin series is not known, it can
be computed via the fast Fourier transform (\acs{FFT}).
Since the Maclaurin series converges, we can find an $L \in \NN$
such that the polynomial $h$ of degree $L$ that we get by truncating the series
after $L+1$ terms approximates $f$ with negligible error.
Thus, the problem simplifies to finding the best approximation
$\tilde{r}^*$ to a polynomial $h$.
The theorem below enables us to compute the best approximation in
$\widetilde{R}_{mn}$ of a polynomial.

\begin{theorem}[Trefethen]
  \label{thm:trefethen}
  The polynomial $h(z) = a_0 + \dots + a_L z^L$ has a unique
  best-approximation $\tilde{r}^*$ out of $\widetilde{R}_{mn}$.
  Let
  \begin{equation*}
    H_{mn} := \begin{pmatrix}
      a_{m-n+1} & a_{m-n+2} & a_{m-n+3} & \cdots & a_L \\
      a_{m-n+2} & a_{m-n+3} & \cdots    & a_{L} \\
      a_{m-n+3} & \cdots    & a_{L} \\
      \vdots    & \iddots \\
      a_{L} &&&& 0
    \end{pmatrix}
    \in \CM{(K + n - m)}{(K + n - m)}
    \,,
  \end{equation*}
  where we define $a_k = 0$ for $k < 0$.
  The error of the approximation $\tilde{r}^*$ is
  \begin{equation*}
    \| h - \tilde{r}^* \|_S = \sigma_{n+1}(H_{mn})
    \,,
  \end{equation*}%
  where $\sigma_{n+1}(H_{mn})$ is the $(n+1)$-st singular value of the
  matrix $H_{mn}$.
  Furthermore, $\tilde{r}^*$ is given by
  \begin{equation}
    \label{eq:cf-error}
    h(z) - \tilde{r}^*(z)
    = \sigma_{n+1} z^L \frac{ u_1 + \dots + u_{K+n-m} z^{L+n-m-1}}
                            { v_{L+n-m} + \dots + v_1 z^{L+n-m-1}}
  \end{equation}
  where $u = (u_1, \dots, u_{L+n-m})^T$ and
  $v = (v_1, \dots, u_{L+n-m})^T$ are the
  $(n+1)$st columns of $U$ and $V$, respectively,
  in the singular value decomposition $H_{mn} = U \Sigma V^H$.
\end{theorem}
\begin{proof}
  See \citep[Theorem~3.2]{Trefethen1981b}.
\end{proof}

This theorem provides us with a formula for the error of the
approximation and we can now work backwards from the error to obtain the
approximation $\tilde{r}^*\in\widetilde{R}_{mn}$ via \eqref{eq:cf-error}.
From $\tilde{r}^*$ we then obtain
the \ac{CF}-approximation $r^\textrm{cf}\in R_{mn}$ by dropping the terms with
negative indices from the numerator.
Since we want to be able to write the rational approximation in
partial fraction decomposition form \eqref{eq:partial-fraction-decomposition},
we restrict ourselves to the case $n = m+1$ in the following.
The whole procedure for this case is listed in Algorithm~\ref{alg:cf}.

\begin{algorithm}
  \caption{Computation of the Carathéodory-Fejér Approximation}
  \label{alg:cf}
  \begin{algorithmic}[1]
    \Function{CF}{$(a_i)_{j=0}^L$, $n$}
      \State $h(z) := \sum_{j=0}^L a_j z^j$
        \Comment{the polynomial to approximate}
      \State
        \label{ln:cf-error-comp-begin}
        $ \displaystyle
          H := \begin{pmatrix}
            a_0 & a_1 & a_2 & \cdots & a_L \\
            a_1 & a_2 & & \iddots \\
            a_2 & & \iddots \\
            \vdots & \iddots \\
            a_L &&&&0
          \end{pmatrix} $
      \State
        $U, \Sigma, V := \Call{SVD}{H}$
        \Comment{$U \Sigma V^H$ is the singular value decomposition of $H$}
      \State
        $\sigma_{n+1} := \Sigma_{n+1,n+1}$
        \Comment{the $(n+1)$-st singular value of $H$}
      \State $u := U_{*,n+1}$ \Comment{the $(n+1)$-st column of $U$}
      \State $v := V_{*,n+1}$ \Comment{the $(n+1)$-st column of $V$}
      \State
        $p(z) := u_1 + u_2 z + \dots + u_{k+1} z^L$
      \State
        \label{ln:cf-error-comp-end}
        $q(z) := v_{K+1} + v_{K} z + \dots + v_1 z^L$
      \State
        \label{ln:cf-r-tilde-star}
        $ \displaystyle
          \tilde{r}^*(z)
            := h(z) - \sigma_{n+1} z^L \frac{p(z)}{q(z)}$
      \State
        \label{ln:cf-denom-comp-begin}
        $z_1, \dots, z_K := \text{the roots of $q$ whose absolute modulus
        is larger than one}$
      \State
        \label{ln:cf-denom-comp-end}
        $q_\textrm{out}(z) := (z - z_1) \cdot (z - z_2) \cdots (z - z_K)$
        \Statex \Comment{the denominator of the \ac{CF}-approximation}
      \State
        $(d_j)_{j=-\infty}^\infty := \Call{LaurentSequence}{q_\textrm{out} \tilde{r}^*}$
        \Comment{$q_\textrm{out}^* \tilde{r} = \sum_{j=-\infty}^\infty d_j z^j$}
      \State
        $\displaystyle
         r^\textrm{cf}(z) := \frac{d_0 + d_1 z + \dots + d_{n-1} z^{n-1}}
                                  {q_\textrm{out}(z)}$
        \Comment{drop the coefficients with negative indices}
      \State
        \Return ($r^\textrm{cf}$, $(z_1, \dots, z_K)$)
    \EndFunction
  \end{algorithmic}
\end{algorithm}

The algorithm starts by computing the quantities used in
Theorem~\ref{thm:trefethen}
(ll.~\ref{ln:cf-error-comp-begin}--\ref{ln:cf-error-comp-end}).
Rearranging the error equation \eqref{eq:cf-error} yields $\tilde{r}^*$,
the best approximation in $\widetilde{R}_{mn}$ (l.~\ref{ln:cf-r-tilde-star}).
Unfortunately, in this form $\tilde{r}^*$ does not provide us with the
coefficients of the numerator and denominator.
Since $h$ is a polynomial, it is easy to see that the poles of $\tilde{r}^*$
are the roots of the polynomial $q$.
We want to write $\tilde{r}^*$ in the form of \eqref{eq:r-tilde-form}.
By definition, all roots of the denominator lie outside of the unit disc,
while the poles inside the unit disc are part of the numerator.
Thus, we obtain the denominator of $\tilde{r}^*$ by multiplying all
linear factors of $q$ corresponding to roots outside the unit disc
(ll.~\ref{ln:cf-denom-comp-begin}--\ref{ln:cf-denom-comp-end}).
Finally, we obtain the coefficients of the numerator by
computing the Laurent series of $q_\textrm{out} \tilde{r}^*$ and then
dropping the terms with negative indices.

Note that $K \le n$ and in general we expect $K$ to be equal to $n$.
In the case where a root of the denominator $q$ lies on the unit circle,
$K$ can be less than $n$. In this case, the root in the denominator is
canceled by a root in the numerator $p$. For all practical purposes,
however, we can assume that by choosing $n$ we can choose the degree $K$ of the
rational approximation \citep{Trefethen1981b}.

\subsubsection{Using the Faber Transform}
\label{sec:faber-cf}

In many practical applications it is desirable to compute approximations
to functions that are defined on domains other than the unit disc.
In our case we are interested in computing an approximation that is
accurate at the eigenvalues of the matrix $M$ \eqref{eq:matrix-exp-solution}.
Since the eigenvalues of $M$ are purely imaginary we can restrict the approximation
domain to an interval on the imaginary axis.
While it would be simple to compute approximations on a disc with a radius
large enough to include the desired interval, being able to choose the
approximation domain more precisely, and hence smaller, often leads to
a better approximation accuracy.

Using the Faber transform, the \ac{CF} approximation can be extended to
allow for the approximation of functions defined on more general domains.
Key to this modification is the observation in \citep{EllFaber1983}
that the Faber transform maps a rational function onto a rational
function. We shall discuss the method introduced in \citep{EllFaber1983}, which we
modify to compute the rational approximation in partial fraction decomposition
form \eqref{eq:partial-fraction-decomposition}.

The Faber transform is based on the fact that Faber polynomials can be used to
derive a series expansion of analytic functions. More precisely, let
$E \subset \mathbb{C}$ be a compact set such that the complement
$E^\textrm{c}$ of $E$ is simply connected in the extended complex plane.
Then, an argument involving the Riemann mapping theorem \citep{GreeneKim2017}
shows that there exists a conformal map $\eta$ that maps the complement
$D^\textrm{c}$ of the closed unit disc conformally onto $E^\textrm{c}$
such that $\eta(\infty) = \infty$ and $\lim_{z\to\infty} \eta(z)/z = c$.
Using this function $\eta$ we can construct a family of polynomials
$p_j$ (for $j = 1, 2, \dots)$ with $p_0(w) = 1$ and $\deg(p_j) = j$,
such that every analytic function $f$ on $E$ can be written as
\begin{equation}
  \label{eq:faber-series}
  f(w) = \sum_{j=0}^\infty a_j p_j(w)
  \quad\text{where}\quad
  a_n = \frac{1}{2\pi \iu} \int_{|z| = 1 + \varepsilon} f(\eta(z)) z^{-n-1}
        \,\df{z}
  \,,
\end{equation}
where $\varepsilon > 0$ has to be chosen small enough
\citep{CurFaber1971,EllFaber1983} such that $g$ is analytic on
$\eta(D_{1+\epsilon})$, where
$D_{1+\epsilon} := \{ z : | z | \le 1 + \epsilon \}$.
These polynomials $p_j$ are called the \emph{Faber polynomials} of $E$. Note that
they only depend on $\eta$ and not on $f$.

Let $g$ be analytic on the unit disc, i.e.,
$g(z) = \sum_{j=0}^\infty c_j z^j$.
The Faber transform of $g$ is given by
\begin{equation*}
  [\mathcal{T} g] (w) = \sum_{j=0}^\infty a_j p_j(w)
  \,.
\end{equation*}%
In other words, the \emph{inverse} Faber transform $\mathcal{T}^{-1} f$ of $f$ is
given by replacing $p_j$ by $z_j$ in the Faber series \eqref{eq:faber-series}
of $f$. Furthermore, the Faber coefficients $a_j$ can be computed without
knowing the Faber polynomials.

As already mentioned, the Faber transform maps rational functions onto
rational functions. Hence, we can obtain an
approximation to a function $f$ defined on $E$ by computing the
\ac{CF} approximation of the inverse Faber transform of $f$ and then
computing the Faber transform of the resulting rational approximation.
The whole method is given in Algorithm~\ref{algo:faber-cf}.

\begin{algorithm}
  \caption{Computation of the Faber-CF-Approximation}
  \label{algo:faber-cf}
  \begin{algorithmic}[1]
    \Function{FaberCF}{$f$, $\eta$, $L$, $n$}
      \State \label{ln:fcf-comp-faber-coeffs}
        $(a_j)_{j=0}^\infty := \Call{FaberSequence}{f}$
        \Comment{$f(w) = \sum_{j = 0}^\infty a_j p_j(w)$}
      \State \label{ln:fcf-compute-cf}
        $(r^\textrm{cf}, (z_1, \dots, z_K)) := \Call{CF}{(a_j)_{j=0}^L, n}$
        \Comment{use Alg.~\ref{alg:cf} here}
      \State $\sigma_j := \eta(z_j)$ \textbf{for} $j = 1, \dots, K$
      \State \label{ln:fcf-transformed-roots}
        $\big(b_j^{(k)}\big)_{j=0}^\infty
          := \Call{FaberSequence}{ w \mapsto \frac{1}{w - \sigma_k}}$
          \textbf{for} $k = 1, \dots, K$
      \Statex
        \Comment{$\sum_{j=0}^\infty b_j^{(k)} p_j(w)
                  = \frac{1}{w - \sigma_k}$}
      \State \label{ln:fcf-linear-system-begin}
        $( c_j )_{j=0}^\infty
          := \Call{Maclaurin}{r^\textrm{cf}} $
        \Comment{$\tilde{r}(z) = \sum_{j=0}^\infty c_j z^j$}
      \State \label{ln:fcf-linear-system-end}
        solve \quad
        $ \left( \begin{array}{@{}c|c|c@{}}
          b_0^{(1)}        & \multirow{4}{*}{\dots} & b_0^{(K)} \\
          b_1^{(1)}        & & b_1^{(K)} \\
          \vdots           & & \vdots \\
          b_{K-1}^{(1)} & & b_{r-1}^{(K)}
          \end{array} \right)
          \begin{pmatrix}
            \beta_1 \\ \beta_2 \\ \vdots \\ \beta_K
          \end{pmatrix}
          =
          \begin{pmatrix}
            c_0 \\ c_1 \\ \vdots \\ c_{K-1}
          \end{pmatrix}$
        \State
          \Return
            $((\sigma_1, \dots, \sigma_K),
              (\beta_1, \dots, \beta_K))$
          \Comment{$ f(w) \approx
                            \sum_{j = 1}^K \frac{\beta_j}{w - \sigma_j}$}
    \EndFunction
  \end{algorithmic}
\end{algorithm}

First, the algorithm computes the coefficient of the Faber series
(l.~\ref{ln:fcf-comp-faber-coeffs}). These are the coefficients of the
Maclaurin series of $\mathcal{T}^{-1} f$.
The algorithm then uses the first $L+1$ coefficients to compute a \ac{CF} approximation
for this analytic function (l.~\ref{ln:fcf-compute-cf}).
It has been shown in \citep{EllFaber1983} that the poles of
$\mathcal{T} r^{\textrm{cf}}$ are $\eta(z_j)$, where $z_j$ for
$j = 1, \dots, K$ are the poles of $r^\textrm{cf}$.
These poles are computed in the next step (l.~\ref{ln:fcf-transformed-roots}).
At this stage of the algorithm we know that the approximation takes the
form $w \mapsto \sum_{j = 1}^K \frac{\beta_j}{w - \sigma_j}$ and it remains to
determine the coefficients $\beta_j$ for $j = 1, \dots, K$.
Since the Faber transform is linear, considering the Maclaurin series of
both sides of the equation
\begin{equation*}
  r^\textrm{cf} =
  \mathcal{T}^{-1} \bigg(w \mapsto \sum_{j = 1}^K \frac{\beta_j}{w - \sigma_j}
                   \bigg)
\end{equation*}
yields a linear system for the coefficient $\beta_j$.
This computation is the final step of the algorithm
(l.~\ref{ln:fcf-linear-system-begin}--\ref{ln:fcf-linear-system-end}).

For the purpose of applying this algorithm to the Schrödinger equation, we need to find
a suitable mapping $\eta$. As mentioned before, the eigenvalues of the matrix $M$ are all
purely imaginary.
Hence, we choose the conformal mapping
\begin{equation}
  \label{eq:map-to-imaginary}
  \eta(z) := \tfrac{R_1}{2} \left(z - \tfrac{1}{z} \right)
  \,,
\end{equation}
which maps the unit sphere $S$ onto the interval $\iu [-R_1, R_1]$
imaginary axis, where all eigenvalues of the matrix $M$ are located.
It is here where the problem at hand needs to be taken into account.
Specifically, if the matrix $M$ has eigenvalues in a different domain, the Riemann mapping $\eta$ needs to be chosen differently.

To summarize, with $\eta$ tailored for the IVP~\eqref{eq:ode-system}, we can
a priori compute shifts $\sigma_j$ and coefficients $\beta_j$, $j=1, \dots, K$,
to compute the rational approximation $r(M)$ in partial fraction decomposition
form \eqref{eq:partial-fraction-decomposition}.
This allows us to approximate the matrix exponential in order to evaluate
the matrix exponential \eqref{eq:matrix-exp-solution} for a given time $\tau$.

\subsection{Step-Size Requirements} 
\label{sec:step-sizes}

In principle, the exponential formula \eqref{eq:matrix-exp-solution} allows
us to compute arbitrary large time steps. There are, however,
practical limitations.
Assume, we choose a large time-step size $\tau$. The solution of the \ac{IVP}
at time $\tau$ is given by $\exp(\tau M)\, u_0$. When $\tau$ is large, the
spectral radius of $\tau M$ is large as well. As discussed in
Section~\ref{sec:rexi} the rational approximation should be close to the true
function values at the eigenvalues of $\tau M$, which makes the
computation of the rational approximation more challenging.

Computing the \ac{CF} approximation for a large domain is more expensive than
computing it for a smaller domain. For a larger domain the degree of the
rational approximation needs to be larger, and the number of terms of the
Maclaurin series that we need also becomes larger, which makes the
computation of the singular value decomposition more and more expensive.
Furthermore, the computation of the \ac{CF} approximation is already expensive
on its own. Hence, we would like to compute the approximation only once and
then reuse it.
This choice, however, limits the step size that our method is able to
perform, as we shall see.

Let us examine the approximation error of the \ac{REXI} method.
For this purpose, let $f: \CN \to \CN$ be a given function. We would like to
compute the matrix function $f(G)$ for some matrix $G$. For an
approximation $r$ to $f$ we define the error function $e$ by $e(\omega) := f(\omega) - r(\omega)$.
Computing $r(G)$ instead of $f(G)$ results in the error
$f(G) - r(G)$ and it is easy to see that $f(G) - r(G) = e(G)$.
Hence, to compute a bound on the approximation error, we need to find a bound
on the norm of $e(G)$. It turns out that if we can bound the error function
$e$ in the eigenvalues of $G$, we can bound $e(G)$.

\begin{proposition}
  \label{pro:error-bound}
  Let $G \in \CM{n}{n}$ and assume that $G = X \diag(\omega_1, \dots, \omega_n) X^{-1}$.
  Furthermore, assume that there exists $\epsilon > 0$ such that
  \begin{equation*}
    e(\omega_j) < \epsilon
    \quad \text{for} \quad
    j = 1, \dots, n.
  \end{equation*}
  Then, $\| e(G) \|_\infty \le \epsilon \cond_\infty(X)$
  where $\| v \|_\infty := \max_{j} | v_j |$, $\| G \|_\infty$ is the
  corresponding operator-norm,
  $\| G \|_\infty = \max_k \sum_{j = 0}^n | G_{kj} |$,
  and
  $\cond_\infty(X) = \| X \|_\infty \, \| X^{-1} \|_\infty$.
\end{proposition}
\begin{proof}
  Follows from \citep[Theorem~4.25]{Higham2008}.
\end{proof}

We can apply this proposition to the case of solving the Schrödinger equation
using the \ac{REXI} method. We know that all eigenvalues of $M$ lie on the
imaginary axis.
Assume that $e(\omega) = f(\omega) - r(\omega) < \epsilon$ for
$\omega \in \iu [-R_1, R_1]$. Then, if we set
\begin{equation}
\label{eq:convergence-condition}
  \tau = \frac{R_1}{\sr(M)}
  \,,
\end{equation}
where $\sr(M) := \max_j | \lambda_j |$ is the spectral radius of $M$,
we have that the eigenvalues of $\tau M$,
which are $\tau \lambda_j$, fulfill that
$| \tau \lambda_j | < R_1$.
Hence, then the assumptions of
Proposition~\ref{pro:error-bound} are satisfied.
Thus, to guarantee proper simulations results we should make sure that
\eqref{eq:convergence-condition} is satisfied.
This restriction is inherent to all variants of the \ac{REXI} methods.

Note that we only need a rough estimate for the largest eigenvalue of $M$
in order to ensure that the accuracy condition
\eqref{eq:convergence-condition} is fulfilled. Such an estimate can be
obtained by running only a few iterations of a sparse eigensolver
\citep[see, e.g.,][and the references
therein]{Stewart_AKrylovSchurAlgorithmForLargeEigenproblems}.
Especially when running many time steps, the time for estimating the largest
eigenvalue can be neglected, since it only needs to be computed once.

Furthermore, note that in order to compute larger time steps, we can choose larger values
for $R_1$ in the conformal map \eqref{eq:map-to-imaginary}.
Choosing larger values of $R_1$, however, is likely to increase the
approximation error $\epsilon$.
Hence, to compensate for an increase in error, one needs to increase the
degree $K$ of the rational approximation, which increases the overall cost
of the method, because it requires more linear systems to be solved.

\begin{table}
  \centering
  \caption{A parameter choice for the Faber-CF-Approximation with
    corresponding properties.}
  \label{tbl:fcf-parameters}
  \begin{tabular}{@{}lllll@{}}
    \toprule
      $K$ & $\eta$ & $R_1$ & Approx.\ Domain & $\| e \|_\infty$ \\
    \midrule
      $16$
        & see eq.~\eqref{eq:map-to-imaginary}
        & $10$
        & $\iu [-10, 10]$
        & \num{2.38e-9} \\
    \bottomrule
  \end{tabular}
\end{table}

In the remainder of this paper, we shall use the Faber-CF-Approximation
defined by the parameters listed in Table~\ref{tbl:fcf-parameters}.
This choice leads to an approximation which has an error of roughly
\num{1e-9} on the approximation interval $\iu [-10, 10]$.
This interval contains about three periods of $\exp(\iu \omega)$. Hence,
using this approximation one time step of the \ac{REXI} method can contain up
to three oscillations of the solution (at a specific point in the spatial
domain).

Comparing this approximation to the approximation derived in
\citep{HautEtAl2016,SchreiberEtAl2018}, we constructed an approximation
using the method from these publications with a comparable accuracy on the
approximation interval $\iu [-10, 10]$. To the best of our knowledge, this
approximation requires at least $34$ poles to achieve the same accuracy.
Thus, while the REXI method using the Faber-CF-Approximation needs to solve
$16$ linear systems per time step, the method from these publications
needs to solve $68$ linear systems for the same time-step size, because
it has to solve two linear systems per pole, as discussed in
Section~\ref{sec:rexi}.

\subsection{Stability}

Since every step of the computation introduces rounding errors and an approximation
error, it is important that these errors are not amplified in the following
steps. Amplification of these errors would cause a run-away effect, which
leads to an exponential growth of the error and needs to be avoided.
Methods that are stable damp the error and thus keep the error under control.

A standard way to asses the \emph{stability} of a method is to apply it to
Dahlquist's test equation
\citep{QuarteroniEtAl_NumericalMethematics,HairerWanner_SolvingOrdinaryDifferentialEquationsII},
which is the \ac{ODE}
\begin{equation}
  \label{eq:test-equation}
  \left\{
  \begin{alignedat}{2}
    \frac{\partial u}{\partial t}(t) &= \omega \cdot u(t)
      & \quad & \text{for $t > 0$} \\
    u(0) &= 1
  \end{alignedat}
  \right.
  \,
\end{equation}
for $\omega \in \CN$.
Let $u_\omega^n$ be the value computed by the method under consideration
after $n$ steps for a particular choice of $\omega$. We define the
\emph{domain of absolute stability} of the method by
\[
  S = \{ \omega \in \CN : | u_\omega^n | \to 0 \quad\text{for}\quad n \to
  \infty \}
  \,.
\]
The process of solving our \ac{ODE} of interest \eqref{eq:ode-system} using a
method with stability domain $S$ and time step $\tau$ is considered stable, if
$\tau \lambda \in S$ for all eigenvalues $\lambda$ of the matrix
$M = (\iu B)^{-1} A$
\citep[cf.][Chapter~IV]{HairerWanner_SolvingOrdinaryDifferentialEquationsII}.
Hence, to analyze the stability of the \ac{REXI} method, we compute the
domain of absolute stability of the method.

To compute the stability domain, we need to solve the test equation
\eqref{eq:test-equation} using the REXI method, which yields the iteration
$u_\omega^{n+1} = r(\omega) \cdot u_\omega^n$, where $r$ is the chosen
rational approximation of the exponential.
Thus, $u_\omega^{n} = (r(\omega))^n \cdot u_\omega^0$, and as a consequence,
$| u_\omega^n | \to 0$ for $n \to \infty$ if and only if $|r(\omega)| < 1$.
Therefore, the stability domain of the \ac{REXI} method is given by
\[
  S = \{ \omega \in \CN : | r(\omega) | < 1 \}
  \,.
\]

\begin{figure}
  \centering
  \subfigure[Stability domain (shaded region) of the \ac{REXI} method.]{%
    \label{fig:fcf-stability-domain}%
    \includegraphics[width=.5\textwidth]{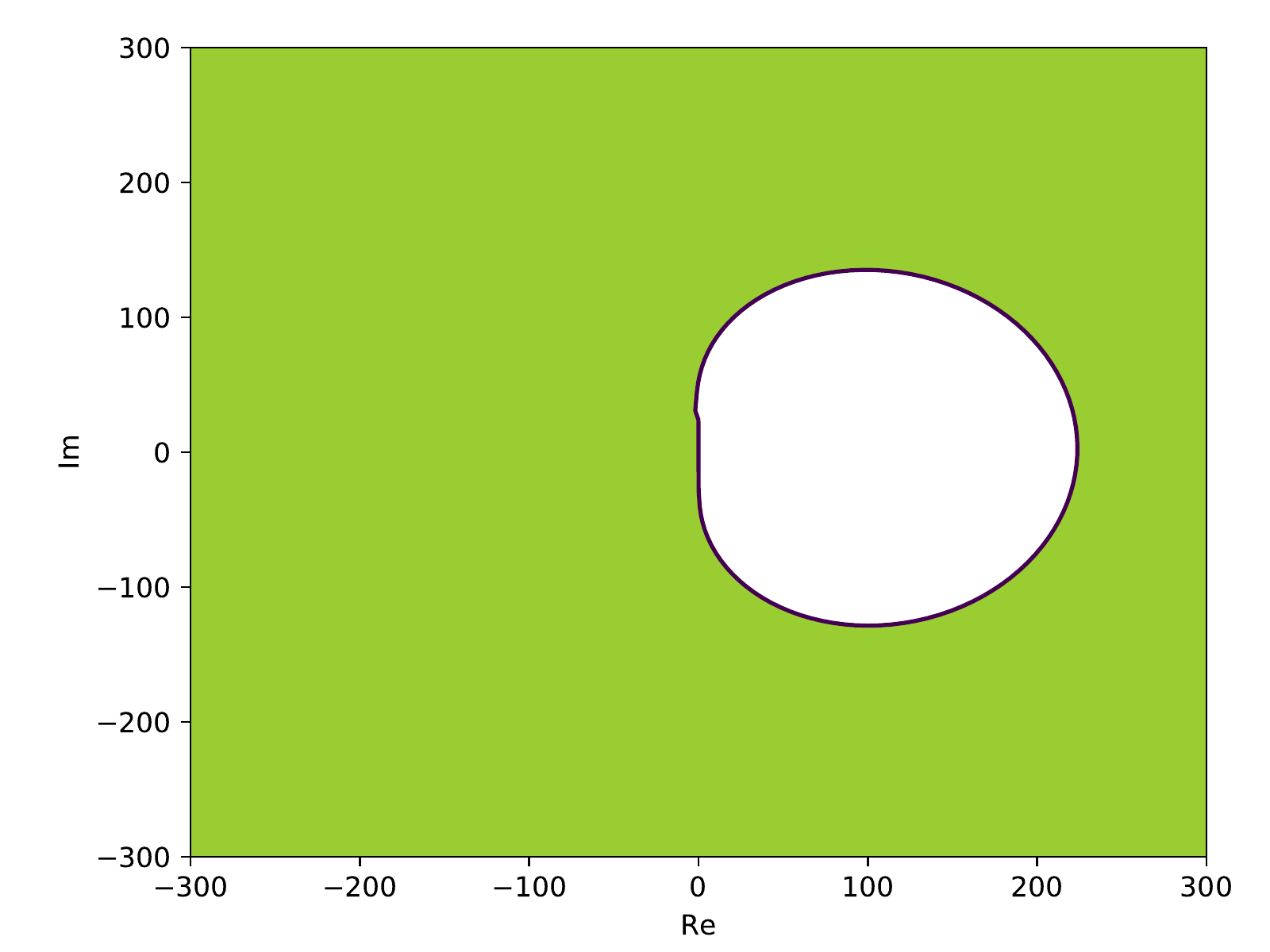}%
  }%
  \subfigure[$| r(z) | - 1$ for $z$ on the imaginary
    axis. The points where these values are smaller than zero belong to
    the stability domain.]{%
    \label{fig:fcf-stability-function-imag}%
    \includegraphics[width=.5\textwidth]{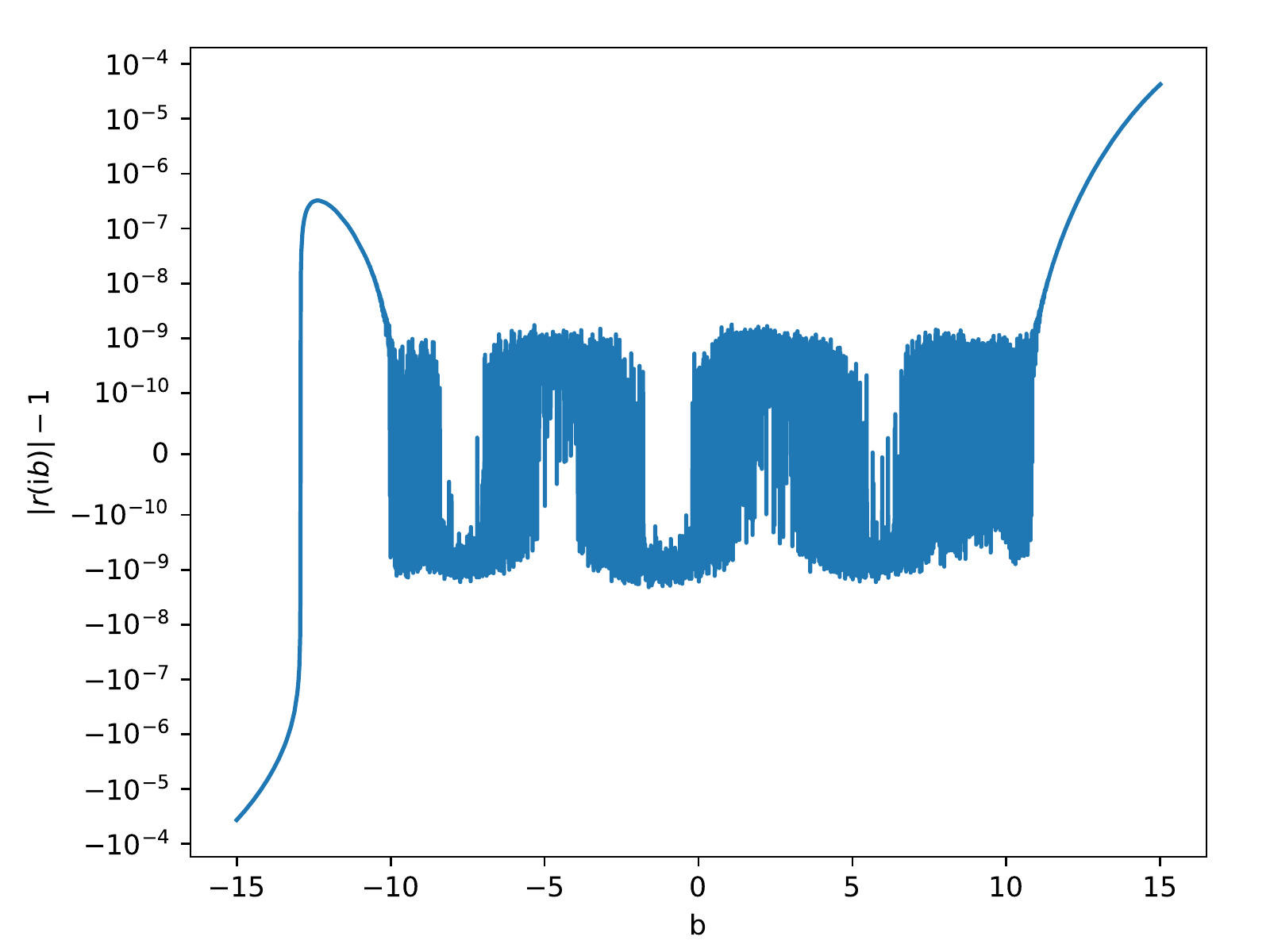}%
  }%
  \caption{%
    Stability of the \ac{REXI} method using the
    Faber-CF approximation specified in Table~\ref{tbl:fcf-parameters}.
  }%
  \label{fig:rexi-stability}
\end{figure}

Turning to the specific \ac{REXI} method considered in this paper, we see in
Figure~\ref{fig:fcf-stability-domain} that the method has a large stability
domain. The eigenvalues of the matrix $M$, however, are located on the
imaginary axis, which is not fully contained in the
stability domain of the method. Taking a closer look at the values of
$| r(z) | - 1$ for $z \in \iu \RN$, as shown in
Figure~\ref{fig:fcf-stability-function-imag}, reveals that $| r(z) |$ exceeds
one by roughly \num{1e-9} on the approximation interval $\iu [-10, 10]$.
That $| r(z) |$ is close to one on the approximation interval is no surprise,
because $r(z)$ approximates $\exp(z)$, and $| \exp(z) | = 1$ on the imaginary
axis. The absolute value of the approximation, however, exceeds one on the
approximation interval by at most $\| e \|_\infty$, the size of the
approximation error. Hence, the method can be stabilized by multiplying
all coefficients $\beta_j$ used in rational approximation
\eqref{eq:partial-fraction-decomposition} by a factor of $1 - \epsilon$,
where $\epsilon$ is slightly larger than $\| e \|_\infty$.
Since this modification introduces an error which is of the same order of
magnitude as the approximation error, the overall accuracy of the method
is only slightly reduced. Recall that in order for all eigenvalues of $M$
to be contained in the approximation interval, the accuracy condition
\eqref{eq:convergence-condition} needs to be fulfilled.

There might be cases where the largest eigenvalues of the matrix $M$ is not
known and estimating its size is too expensive. Since, the degree of the
numerator is smaller than the degree of the denominator of the
rational approximation $r$, $| r(z) | \to 0$ for $|z| \to \infty$ and thus
large eigenvalues have a good chance of falling into the region of stability.
Hence, if the corresponding mode is irrelevant for the solution process,
usable results might still by obtained. Nevertheless, we recommend to compute
a rough estimate for the largest eigenvalue and use the accuracy condition
\eqref{eq:convergence-condition} to choose the step size.
Note, this behavior is in contrast to polynomial approximations,
like the Chebyshev method \citep{TalEzerKosloff1984}, we will discuss below.
For a polynomial (of degree larger than zero), $| p(z) | \to \infty$ for
$| z | \to \infty$, and hence the method becomes almost immediately unstable
if one of the eigenvalues of $M$ lies outside of the approximation region.

Turning to the specific approximation as specified by
Table~\ref{tbl:fcf-parameters}, we observe
that $| r(z) |$ exceeds one by only $10^{-9}$, which
results in an error growth proportional to
$(1 + 10^{-9})^n
  = e^{n \cdot \log(1 + 10^{-9})}
  \approx e^{n \cdot 10^{-9}}$.
Even though the error grows exponentially, since the base is only slightly
larger than one, it requires a large number of time steps for the
exponential effect to dominate. For example, for an error amplification of
a factor of $2$ we need to run the method for \num{7e8} time steps.
Here, we do not plan to run our method for such a large number of time steps.
Hence, in the following, we do not stabilize the method for the benefit of a
higher accuracy.


\section{Numerical Experiments}
\label{sec:num}

We carry out numerical experiments to study the potential,
effectiveness and efficiency of the method presented above.
To test the algorithms in a realistic setting, we simulate two different
quantum mechanical systems that feature two famous quantum mechanical
phenomena that cannot be explained by classical mechanics.

To this end, we use an implementation of the \ac{REXI} method written
in C++ and parallelized using \ac{MPI}. For the finite element discretization
we use libMesh \citep{KirkEtAl2006} and all matrix and vector operations
are implemented using the PETSc library \citep{Petsc}.
All computations are performed on the \ac{JURECA} cluster \citep{JSC2018},
which consists of 1872 nodes connected via InfiniBand. The nodes we use
contain a two-socket board equipped with two Intel Xeon E5-2680 v3 and
$128 \mathrm{GiB}$ of memory.

In many applications we want to start the
simulation with a particle at a certain location and with a certain momentum
and see how it evolves in time. Due to the Heisenberg uncertainty relation,
however, we can either give a quantum particle a defined position or
a defined momentum, but not both. Hence we choose a Gaussian-wave package as
initial condition.

A Gaussian wave-package is defined as
\begin{equation}
  \label{eq:wave-package}
  \psi(\vec{r}) := C_\mathrm{norm} \cdot
  e^{-\frac{1}{2} (\vec{r} - \bar{\vec{r}})^T \Sigma^{-1} (\vec{r} - \bar{\vec{r}})}
  \cdot e^{\iu \dotp{\bar{\vec{p}}}{\vec{r} - \bar{\vec{r}}} / \hbar}
\end{equation}
where $\Sigma \in \RM{d}{d}$ is a symmetric, positive definite matrix and
$C_\mathrm{norm} \in \mathbb{R}$ is chosen such that
\begin{equation*}
  \int_\Omega \overline{\psi} \psi \df{\vec{r}} = 1
  \,.
\end{equation*}
The wave-package describes a particle ensemble with position expectation
value $\bar{\vec{r}}$ and momentum expectation value $\bar{\vec{p}}$. The matrix
$\Sigma$ describes the uncertainty of the particle position. At the same
time the matrix influences the uncertainty of the momentum---the smaller
the uncertainty of the position the larger the uncertainty of the
momentum and vice versa.

For simplicity, all quantities are measured in Hartree atomic units
\citep{Hartree1928a,MillsEtAl1993}, i.e., we choose a system of measurement
in which the electron mass $m_e$, the elementary charge $e$, the
reduced Plank constant $\hbar$, and the inverse Coulomb constant
$4 \pi \epsilon_0$ are all equal to one. In this system, length is measured
in \emph{bohr}, $a_0$, i.e., the Bohr radius, and energy is measured in
Hartree, $E_h = \hbar / (m_e a_0^2)$.

\subsection{Quantum Tunneling} 

We start our numerical investigation by considering a simulation of
quantum particle tunneling through a step-potential barrier.
Quantum tunneling describes a phenomenon in which a quantum particle passes
through a potential barrier even though its kinetic energy is smaller than
the height of the barrier
\citep[see, e.g.,][]{Messiah1999,Shankar2014,Griffiths2017}.
This behavior is in contrast to classical mechanics where such a behavior
is not possible.

\begin{table}
  \caption{Parameters of the quantum tunneling simulation.}
  \label{tbl:tunneling-parameters}
  \centering
  \begin{tabular}{@{}llllllll@{}}
    \toprule
    Parameter
    & $m$
    & $V_\mathrm{max}$
    & $C_\mathrm{barr}$
    & $\Sigma$
    & $\sigma$
    & $\bar{\vec{r}}$
    & $\bar{\vec{p}}$
    \\ \midrule
    Value
    & $1\,m_e$
    & $15\,E_h$
    & $0.005\, a_0$
    & $\sigma I$
    & $4\, a_0^2$
    & $-3\, a_0$
    & $5 \, \tfrac{\hbar}{a_0}$
    \\ \bottomrule
  \end{tabular}
\end{table}

We consider the tunneling process defined as follows.
An electron moves in the step-potential given by
\begin{equation*}
  V(\vec{r}) = \begin{cases}
    V_\mathrm{max} & \text{if } | r_1 | \le \frac{C_\mathrm{barr}}{2} \\
    0 & \text{otherwise}
  \end{cases}
  \,,
\end{equation*}
which has a barrier at $x = 0$. 
The electron starts at the left of the barrier and moves to the right,
with a speed typical for electrons that are emitted by electron guns.
When the electron reaches the barrier it has a certain probability
of being reflected from the barrier or tunneling though the barrier.
We simulate this process by numerically solving the Schrödinger equation
\eqref{eq:schroedinger} with the step-potential, a Gaussian wave-package
\eqref{eq:wave-package} as initial values, and the parameters listed in
Table~\ref{tbl:tunneling-parameters}.

\begin{figure}
  \centering
  \subfigure[$t = 0\, \hbar/E_h$]{%
    \label{fig:tunneling-t0}
    \includegraphics[width=0.47\textwidth]{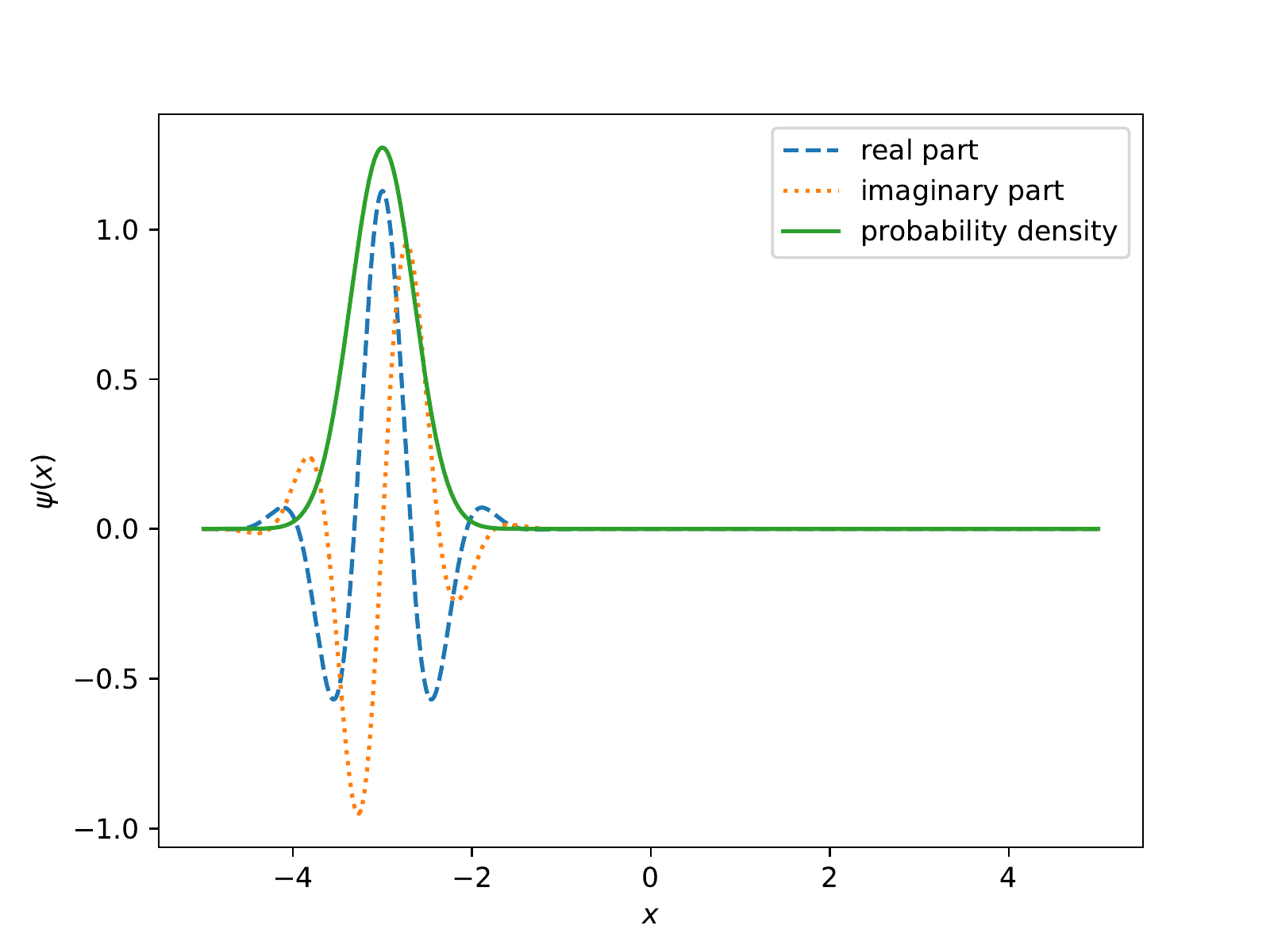}}
  \hfill
  \subfigure[$t = 0.012\, \hbar/E_h$]{%
    \label{fig:tunneling-t1}
    \includegraphics[width=0.47\textwidth]{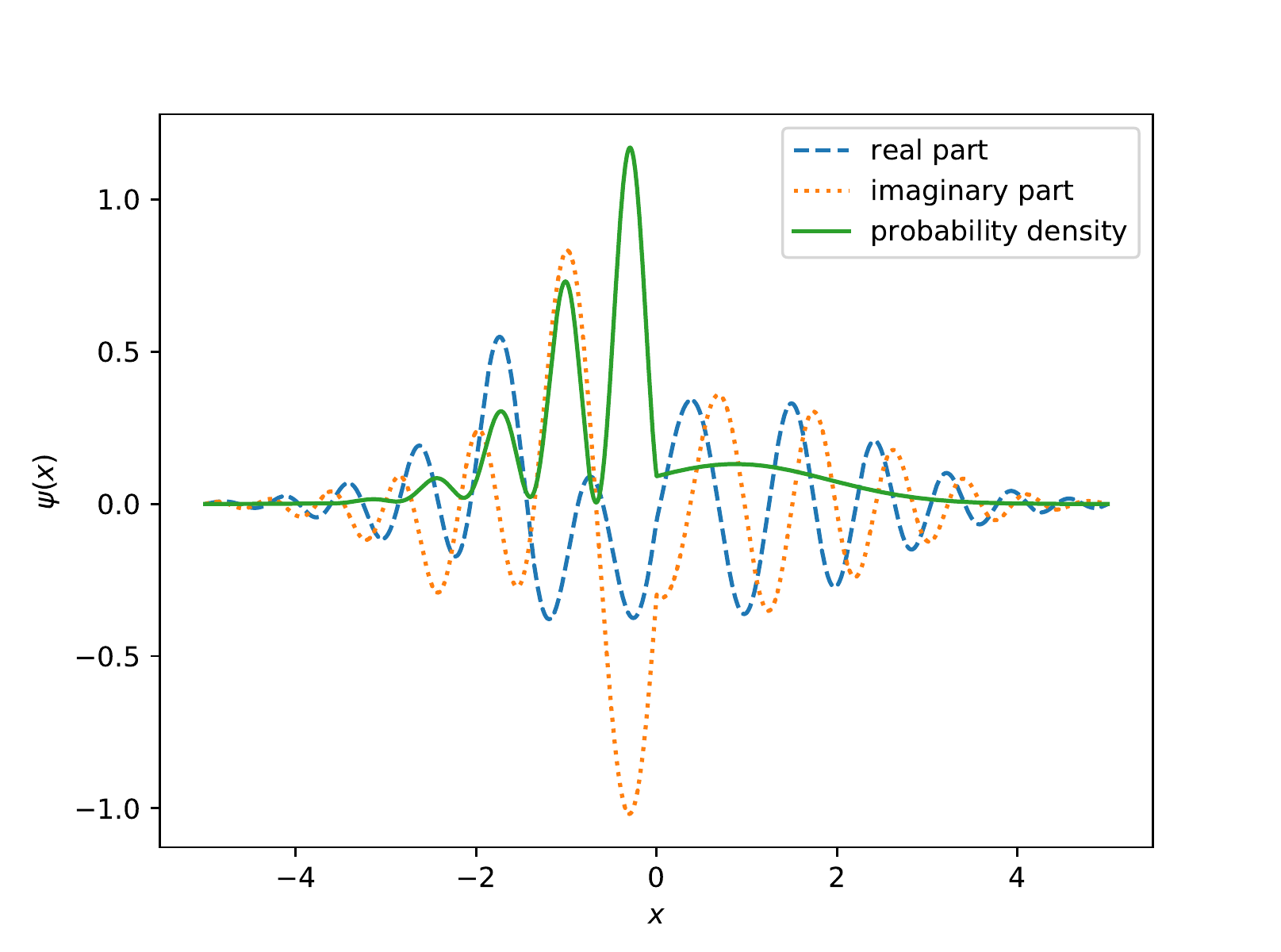}}
  \caption{A 1D tunneling problem. Simulation of the Schrödinger equation
    using the parameters from Table~\ref{tbl:tunneling-parameters}
    and a Gaussian wave-package as initial values at two different times.}
  \label{fig:tunneling}
\end{figure}

Figure~\ref{fig:tunneling} shows the state of the quantum system at two
different times. In the initial state (Figure~\ref{fig:tunneling-t0}) the
probability density is concentrated at the left of the domain.
At a later time (Figure~\ref{fig:tunneling-t1}) a part of the particle
ensemble has passed through the barrier at $x=0$, which can be seen by a raise
of the probability density for $x > 0$. A large portion of the particles,
however, is reflected from the barrier, which results in the interference
pattern caused be the superposition of the incoming and reflected waves.

For the numerical simulation, we discretize the equation using finite
elements as described in Section~\ref{sec:fem}. We choose $4000$ equally sized
finite elements of order two, which results in an \ac{ODE} system with $8001$
degrees of freedom. We then simulate this system using the \ac{REXI}
time-stepping scheme (Section~\ref{sec:rexi}).

In our first experiment, we compare the serial execution time and accuracy
of the \ac{REXI} method to other \ac{ODE} solvers.
This comparison is important, because we later want to investigate the
parallel performance of the \ac{REXI} method, and thus we need to know the
fastest serial method as a reference point, to get a realistic
impression of the effectiveness of the parallelization.

Using \ac{REXI} requires the construction of a rational approximation of the
exponential function. We compute the Faber-CF approximation
(Section~\ref{sec:faber-cf}) as
defined by the parameters listed in Table~\ref{tbl:fcf-parameters} and
discussed earlier.

The first method that we compare the \ac{REXI} method with is the Chebshev
method \citep{TalEzerKosloff1984}.
It works by approximating $\exp(z)$ via
$p(z) := \sum_{k=0}^N a_k T_k(-\iu z/R)$ on the interval $-\iu[-R, R]$, where
$T_k$ ($k \in \NN$) is the Chebyshev polynomial of degree $k$
\citep{Rivlin1981,Powell1981} and $a_k \in \RN$.
The coefficients $a_k$ can be efficiently computed using the \ac{FFT}
\citep{Trefethen2012}.
The polynomial is then used to approximate $\exp(\tau B^{-1}A) u_0$ by
evaluating $p(\tau B^{-1} A) u_0$,
which can be done using the Clenshaw algorithm \citep{Clenshaw1955}.
In contrast to the \ac{REXI} method, the Chebyshev method only allows for
spatial parallelization.
Note that because of the mass matrix $B$, stemming from the finite elements
approach, the application of the Chebyshev method also involves the solution
of linear systems, making it significantly more expensive than in the case of
$B=I$.
For the sake of a meaningful comparison, we match the approximation quality of
the Chebyshev polynomial with the one of the Faber-CF approximation.
We choose the approximation interval $\iu [-10, 10]$ and
a polynomial of degree $26$, which leads to an error of
the same order of magnitude as the rational approximation.

The second method that we compare the \ac{REXI} method with is a
fourth order Rosenbrock method. More precisely, we choose the
$L$-stable method listed in
\citep[Section IV.7, Table~7.2]{HairerWanner_SolvingOrdinaryDifferentialEquationsII}.
Rosenbrock methods are diagonally implicit
Runge-Kutta methods, and hence the method requires the solution of
four linear equations per time step.
We use the implementation of this method provided in PETSc \citep{Petsc}.

There are further methods for numerically simulating the \ac{ODE} system
arising from the discretization of the Schrödinger equation, e.g.,
the Crank-Nicolson \citep{McCulloughWyatt1971a} or
the leapfrog \citep{AskarCakmak1978} method
(see also \citep{LeforestierEtAl1991}).
These methods, due to their low order, however, require very small time steps.
Hence, we do not consider them in the comparison.

\begin{table}
  \centering
  \caption{%
    Tunneling problem timings, $LU$ decomposition in space.
    Comparison of different time-stepping schemes in serial.
    The time in parenthesis is the time needed for solving
    linear systems involving the mass matrix $B$.
    }
  \label{eq:tunneling-method-comparison}
  \begin{tabular}{@{}lllll@{}}
    \toprule
    & $\Delta t$ & error & time/\si{s} \\
    \midrule
    Chebyshev & \num{2e-4} & \num{3.17e-06}   & \PZZ 7.49 (3.08) \\[.4em]
    REXI      & \num{2e-4} & \num{6.66e-07} & \PZZ 2.59 \\[.4em]
    Rosenbrock 4 & \num{2e-4} & \num{1.11e-05} & \PZZ 4.50 \\
              & \num{2e-5} & \num{2.43e-06} & \PZ 44.80 \\
              & \num{5e-6} & \num{4.23e-07} & 179.56 \\
    \bottomrule
  \end{tabular}
\end{table}

We simulate the $1$D tunneling problem using the three different methods
for a time period of $0.2 \hbar/E_h$, and measure the time the different
methods require for the time stepping.
All linear systems are solved using the $LU$ decomposition,
implemented in the SuperLU\_DIST software package \citep{LiDemmel2003},
and all computations are performed sequentially.
The results are listed in Table~\ref{eq:tunneling-method-comparison}.
Note that we choose $\Delta t$ for the REXI and Chebshev method such that
the accuracy condition \eqref{eq:convergence-condition} is fulfilled.

Considering the results, we see that the \ac{REXI} method is the overall
fastest method for this simulation. Furthermore, when taking accuracy into
account, the Rosenbrock method needs substantially more time steps to reach
the same accuracy as the \ac{REXI} method.
Note that the Rosenbrock method has to solve only four linear systems, while
the \ac{REXI} method has to solve 16. Hence, we would expect that the
Rosenbrock method would be four times faster per time step than the \ac{REXI}
method, while we measure it to be two times slower.
We assume that the Rosenbrock method shows this behavior, because we used the
generic implementation provided by PETSc, and a more specialized
implementation would perform better. Nevertheless, due to the smaller
step-size requirement of the Rosenbrock method, the \ac{REXI} method would
still be the fastest method.
Justified by these results, we shall use the time of the serial execution
of the \ac{REXI} method as the reference point for computing the
parallel speedup of the method.

In the second experiment, we want to investigate the parallelization potential
of the \ac{REXI} method. In this experiment we simulate the same quantum
system. To obtain more meaningful time-measurements and to give the
method enough work that can be distributed along multiple processors,
we increase the number of degrees of freedom by refining the mesh four times.
Each refinement split one mesh cell into two and thus roughly doubles the
degrees of freedom each time. Bearing in mind the accuracy condition
\eqref{eq:convergence-condition}, we have to reduce the step size of the
\ac{REXI} method. We use a time step size of $\num{5e-7} \hbar/E_h$ and
simulate the system for $\num{5e-4} \hbar/E_h$.
Recall that we have two types of parallelization that we can use---time
and space-parallelization. For now we restrict ourselves to inspecting the
time-parallelization only.

Note that the time-parallelization is limited by the number of poles of the
rational approximation, because the number of poles determine the maximum number $K$ of
linear systems in \eqref{eq:rational-evaluation} that can/need to be solved
simultaneously. Hence, using the approximation described above allows us
to split the computation of one time step into $16$ independent tasks.

\begin{table}
  \centering
  \caption{%
    Tunneling problem timings, $LU$ decomposition in space, four refinements.
  }
  \label{tbl:tunneling-rexi-scaling}
  \begin{tabular}{@{}rrrrrrrr@{}}
  \toprule
   && \multicolumn{4}{c}{time/\si{s}} \\
   \cmidrule(lr){3-6}
   nodes &  cores & total & rhs & local & reduce & speedup & efficiency \\
  \midrule
     1 &      1 &  39.92 &   1.38 &   38.51 &  0.01 &  1.00 &  1.00 \\[.4em]
     1 &      2 &  26.12 &   1.49 &   23.87 &  0.76 &  1.53 &  0.76 \\
     1 &      4 &  19.20 &   1.70 &   16.29 &  1.21 &  2.08 &  0.52 \\
     1 &      8 &  16.21 &   2.74 &   12.01 &  1.45 &  2.46 &  0.31 \\
     1 &     16 &  15.06 &   3.33 &    7.82 &  3.90 &  2.65 &  0.17 \\[.4em]
     2 &      2 &  21.64 &   1.36 &   19.52 &  0.75 &  1.84 &  0.92 \\
     4 &      4 &  12.11 &   1.31 &    9.73 &  1.07 &  3.30 &  0.82 \\
     8 &      8 &   8.13 &   1.31 &    4.92 &  1.89 &  4.91 &  0.61 \\
    16 &     16 &   5.68 &   1.27 &    2.37 &  2.03 &  7.03 &  0.44 \\
  \bottomrule
  \end{tabular}
\end{table}

The results of the experiment are listed in
Table~\ref{tbl:tunneling-rexi-scaling},
which contains the time measurements for running the method with different
numbers of nodes. In the first half of the table, we keep the number of
compute nodes constant, which means that all tasks are running on the same
two-socket system and thus have direct access to the same memory.
In the second half of
the table, we use one compute node per tasks, which means that each tasks
has its own \ac{CPU} and memory.
In addition to the total runtime, the times for the individual phases of the
algorithm are also given. The \ac{REXI} method evaluates
\eqref{eq:rational-evaluation} in three phases. First, the method
has to compute the right-hand side (rhs) of the linear systems by
multiplying the matrix $\iu B$ and the vector $u_0$. Note that from the
time-parallelization perspective, this is a sequential part of the
algorithm. Second, each process performes the local part of the computation,
i.e., it solves the local linear systems and the local sums.
Third, all processes compute the the global sum and distribute the result
to all processes (reduce), which is the step that involves communication.
In addition to the time measurements, the table contains the speedup and
the parallel efficiency.

Inspecting the table, we see that the efficiency when running on one node is
low. The time for the computation of the right-hand side increases
with increasing number of cores. The time for the local computation achieves
only a speedup of $5$ when running on $16$ cores. Furthermore, the time
needed for the reduction increases. This behavior is due to the fact that
modern \acp{CPU} can run at higher clock speeds when only a few cores are
used, and that all cores on one node share the same memory interface, which
becomes a limiting factor.

Using multiple nodes, the method scales better. We observe
that the time for the right-hand side computation remains constant, which
is expected, because it is the sequential part of the algorithm.
The time for the local operations scales perfectly. The time for the global
summation, however, increases.
Hence, due to the sequential part and the increased communication cost, we
only achieve an efficiency of $0.44$.
Note however, that this algorithm is not meant to be used alone. It should
be used to provide additional parallelism in the situation where increasing
spatial parallelism is not feasible anymore.
With respect to these considerations the speedup is promising.
As a next step, we considered a larger problem and combine temporal and
spatial parallelism.

\subsection{Double Slit Experiment} 

For the purpose of applying the \ac{REXI} method to a larger problem, we
consider the simulation of a double-slit interference experiment with
electrons.
\begin{figure}
  \centering
  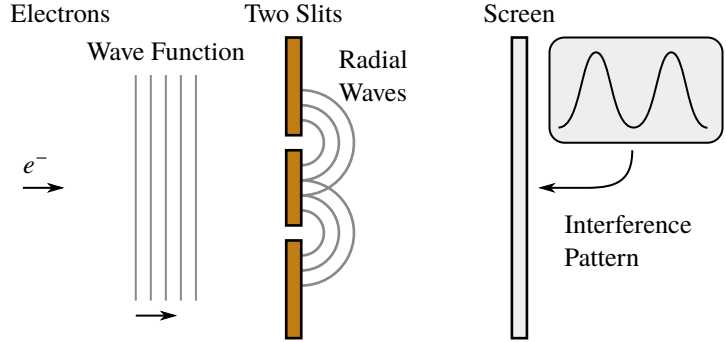
  \caption{Double-slit experiment.}
  \label{fig:double-slit-experiment}
\end{figure}
This experiment demonstrates the
wave-like character of matter particles and shows the limits of classical
mechanics \citep{Shankar2014}.
Assume that we are shooting electrons at a wall that has two thin slits, each of
which can be closed.  Most electrons will hit the wall, some, however, will
pass through one of the slits. If we place a fluorescent screen at the other
side of the wall, we can record the probability density of the incidenting
electrons.  We repeat this experiment three times.  Once with both slits open,
once where the first slit is closed, and once where the second slit is closed.

Classical mechanics would predict that the probability density we measure
with both slits open is just the sum of the probability density that we obtain
in the two cases where just one slit is open. It turns out, however, that we
observe an interference pattern in the case of two open slits.
This interference pattern can be explained using quantum mechanics. We
describe the incidenting electrons by a planar probability wave that moves
in the direction of the screen. When the planar wave hits the wall, each of
the two slits emits radially outward going waves. These waves interfere,
and when the electrons hit the screen, the probability density that results
from this interference becomes visible.
Figure~\ref{fig:double-slit-experiment} shows a schematic overview of the experiment.

Note that this experiment has never been actually performed in precisely this way.
It resembles, however, the essential features of many experiments that have
been performed without the technical complications they involve.
Yet, Tonomura et al.\ conducted an experiment very
close\footnote{%
  Instead of a wall with two slits an electron biprism was used.}
to the one that we described~\citep{TonomuraEtAl1989} and that we shall simulate.

\begin{table}
  \centering
  \caption{Parameters of the double-slit experiment simulation.}
  \label{tbl:double-slit-parameters}
  \begin{tabular}{@{}llllll@{}}
    \toprule
    Parameter
    & $m$
    & $\Sigma$
    & $\sigma$
    & $\bar{\vec{r}}$
    & $\bar{\vec{p}}$
    \\
    \midrule
    Value
    & $1\,m_e$
    & $\sigma I$
    & $10^{-4}\, a_0^2$
    & $(0, -350)\,a_0$
    & $(0, 0.1)\,\tfrac{\hbar}{a_0}$
    \\
    \bottomrule
  \end{tabular}
\end{table}

\begin{figure}
  \centering
  \subfigure[$\Re \psi(\vec{r},t)$ at $t = 0\,\hbar/E_h$]{%
    \includegraphics[width=.5\textwidth]{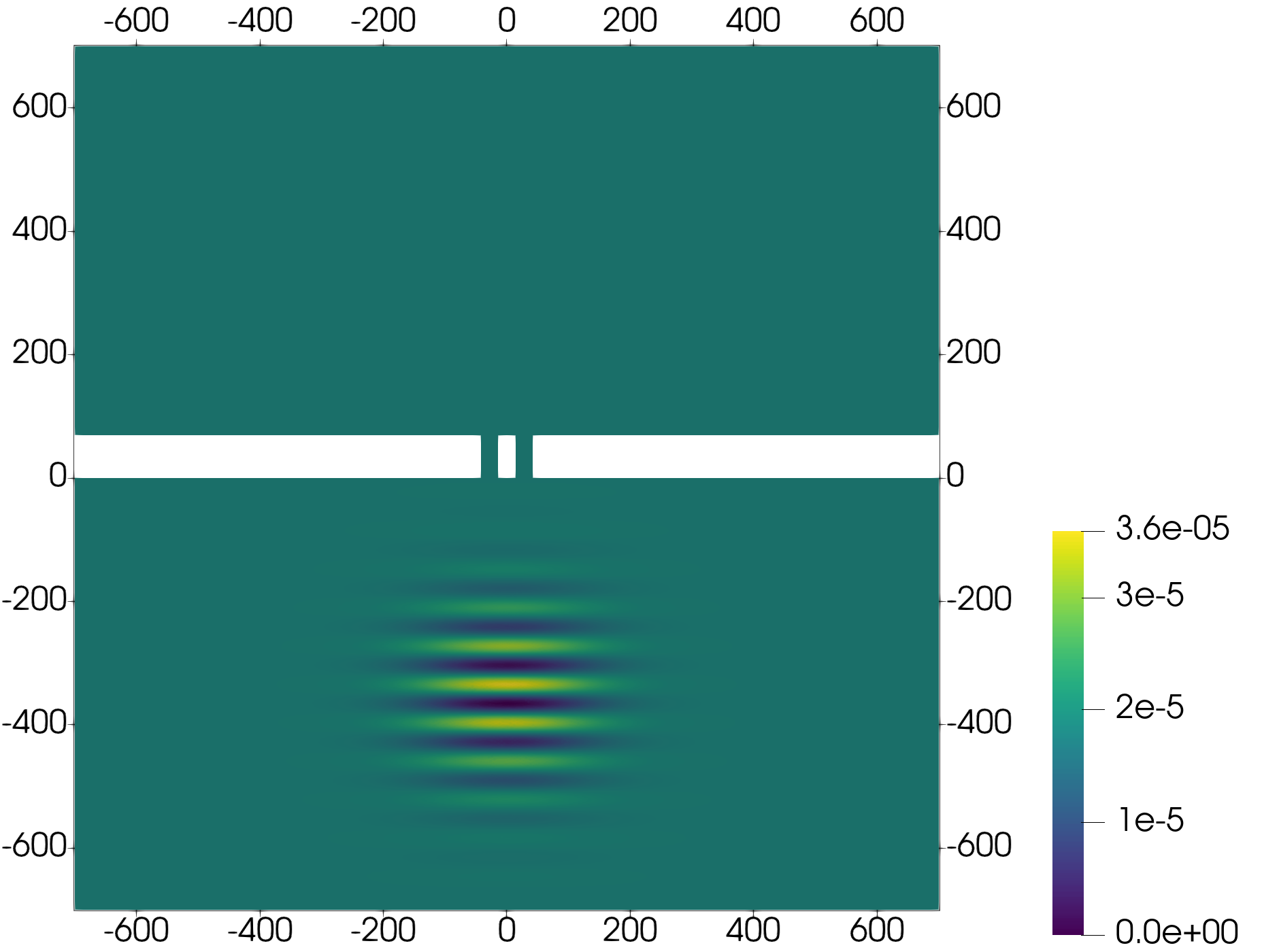}%
  }%
  \subfigure[$\overline{\psi(\vec{r},t)}\,\psi(\vec{r},t)$ at $t = 0\,\hbar/E_h$]{%
    \includegraphics[width=.5\textwidth]{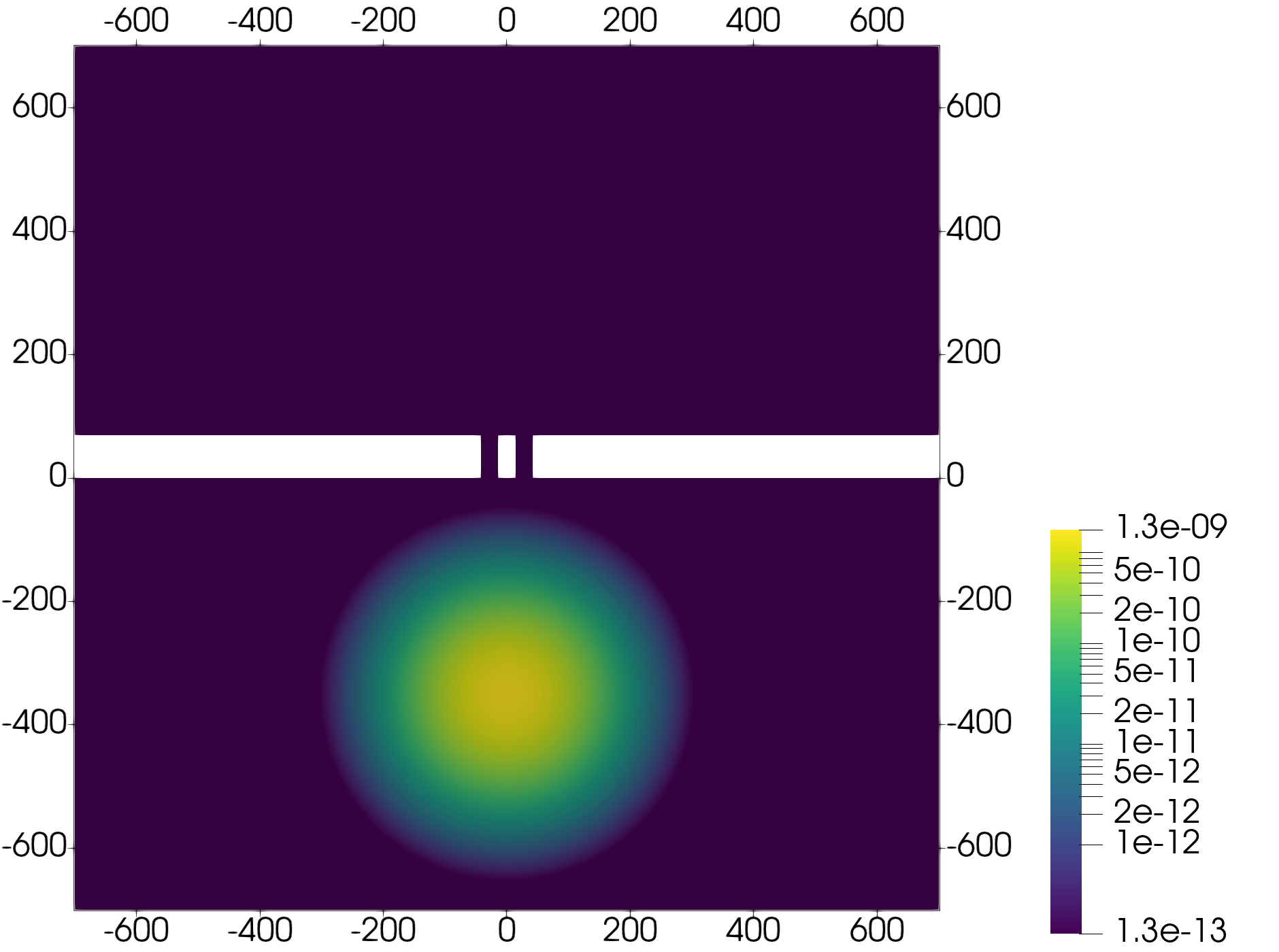}%
  }%
  \\
  \subfigure[$\Re \psi(\vec{r},t)$ at $t = \num{6e5}\,\hbar/E_h$]{%
    \includegraphics[width=.5\textwidth]{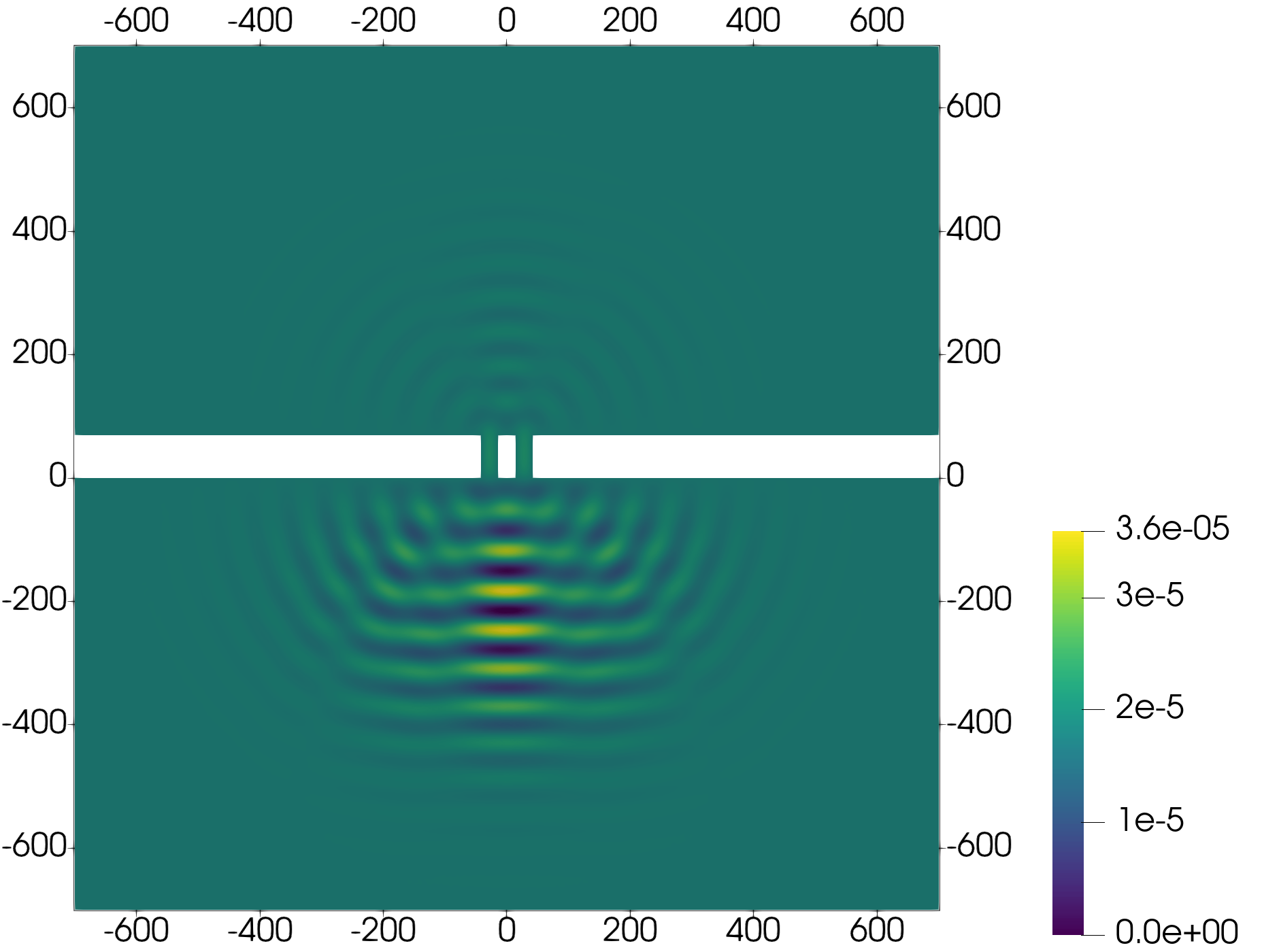}%
  }%
  \subfigure[$\overline{\psi(\vec{r},t)}\,\psi(\vec{r},t)$ at $t = \num{6e5}\,\hbar/E_h$]{%
    \includegraphics[width=.5\textwidth]{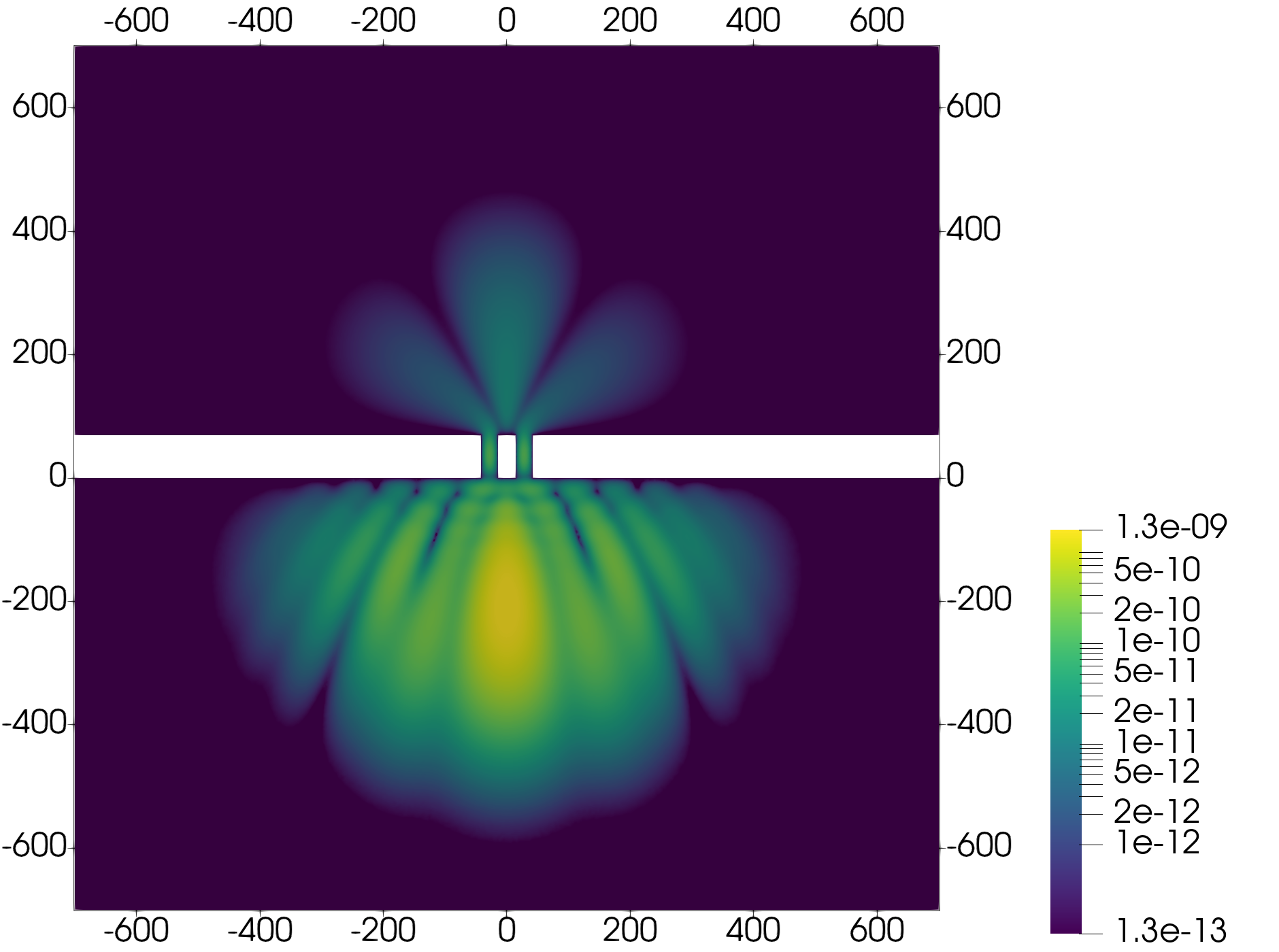}%
  }%
  \caption{%
    Simulation of the Schrödinger equation using parameters from
    Table~\ref{tbl:double-slit-parameters} at two different times.
    The real part (left) and the probability density function (right) are
    shown.
  }
  \label{fig:ds-wave-functions}
\end{figure}

For the simulation, we need to determine the appropriate parameters of the
Schrö\-dinger equation.  We can model the wall with the two slits by choosing
the domain of the \ac{PDE}
appropriately (see Figure~\ref{fig:ds-wave-functions}), imposing zero Dirichlet boundary conditions. These conditions imply that the
particle has a zero possibility of reaching the boundary of the domain and,
hence, must be contained within. Since the electron is supposed to move freely
within the domain we choose the zero potential, $V(\vec{r}) = 0$.
Furthermore, to fulfill the condition \eqref{eq:convergence-condition} we
chose a step size of $10\, a_0$. The remaining
parameters are listed in Table~\ref{tbl:double-slit-parameters}.

The state of the simulation is shown in Figure~\ref{fig:ds-wave-functions} at
two different times. At the beginning the wave package is located in the lower
part of the domain and moving towards the wall. At the second time most
of the particles have been reflected from the wall, but a fraction of them
have passed through either of the slits. At the other side of the wall an
interference pattern forms.

We discretize the equation using the finite element method as described in
Section~\ref{sec:fem}. We use
a suitable triangulation
to obtain a discretization using $179620$ finite elements of order two,
leading to an \ac{ODE} system with $288796$ unknowns. Note that in this case
the finite element method is by far the preferred discretization due to the
complicated geometry of the domain.

\begin{table}
  \centering
  \caption{REXI with LU}
  \label{tbl:rexi-ds-lu-timings}
  \begin{tabular}{@{}r@{$\,\times\,$}lrrrr@{}}
  \toprule
    \multicolumn{2}{c}{cores} &
    \multicolumn{2}{c}{time / \si{s}} \\
    \cmidrule(r){1-2} \cmidrule(lr){3-4}
    time & space    &   total   & reduce &   speedup & efficiency \\
  \midrule
       1 &        1 &    132.85 &  0.00 &    1.00 &       1.00 \\
       1 &        2 &    123.77 &  0.00 &    1.07 &       0.54 \\
       1 &        4 &    114.43 &  0.00 &    1.16 &       0.29 \\
       1 &        8 &     99.70 &  0.00 &    1.33 &       0.17 \\
       1 &       12 &     94.68 &  0.00 &    1.40 &       0.12 \\
       1 &       16 &     97.45 &  0.00 &    1.36 &       0.09 \\
       2 &        1 &     78.00 &  0.04 &    1.70 &       0.85 \\
       4 &        1 &     46.74 &  0.07 &    2.84 &       0.71 \\
       8 &        1 &     27.14 &  0.20 &    4.90 &       0.61 \\
      16 &        1 &     17.42 &  1.51 &    7.63 &       0.48 \\
  \bottomrule
  \end{tabular}
\end{table}

Our first numerical experiment for this setting aims at comparing the
space- and the time-parallelization. We use, again, SuperLU\_DIST to solve
all linear systems in the \ac{REXI} method. The results can be found in
Table~\ref{tbl:rexi-ds-lu-timings}.
We see that the space-parallelization achieves a speedup of barely $1.36$ when using
$16$ cores.
This is due to the fact that SuperLU\_DIST does not seem to parallelize well for this problem size, but we did not investigate this issue further.
Furthermore, we can see that the time-parallelization is much
more effective than the space-parallelization: using $16$ cores in time gave a speedup of about $7.63$, about $5$ times as high as before.
To avoid the complications induced by the LU solver and to compare to a more practical choice, for the next numerical experiment we replace the LU solver by the \acs{GMRES} method \citep{SaadSchultz1986}, which is an iterative solver.

\begin{table}
  \centering
  \caption{Double Slit, Solving with GMRES}
  \label{tbl:rexi-ds-gmres-timings}
  \begin{tabular}{@{}r@{$\,\times\,$}lrrrrr@{}}
  \toprule
    \multicolumn{3}{@{}c@{}}{cores} &
    \multicolumn{2}{c}{time / \si{s}} \\
    \cmidrule(r){1-3} \cmidrule(lr){4-5}
    time & space & total &        total &  reduce & speedup & efficiency \\
  \midrule
       1 &        1 &      1 &    299.74 &   0.00 &      1.00 &       1.00 \\
       1 &       12 &     12 &     46.62 &   0.00 &      6.43 &       0.54 \\
       1 &       24 &     24 &     21.29 &   0.00 &     14.08 &       0.59 \\
       1 &       48 &     48 &      7.81 &   0.00 &     38.40 &       0.80 \\
       1 &       96 &     96 &      2.97 &   0.00 &    100.87 &       1.05 \\
       1 &      192 &    192 &      2.01 &   0.00 &    149.08 &       0.78 \\
       1 &      384 &    384 &      1.98 &   0.00 &    151.15 &       0.39 \\
       8 &       24 &    192 &      2.84 &   0.09 &    105.72 &       0.55 \\
       8 &       48 &    384 &      1.18 &   0.02 &    254.50 &       0.66 \\
       8 &       96 &    768 &      1.05 &   0.04 &    286.56 &       0.37 \\
      16 &       12 &    192 &      3.12 &   0.13 &     96.04 &       0.50 \\
      16 &       24 &    384 &      1.57 &   0.11 &    190.62 &       0.50 \\
      16 &       48 &    768 &      0.89 &   0.08 &    335.15 &       0.44 \\
  \bottomrule
  \end{tabular}
\end{table}

The results of the second numerical experiment can be found in
Table~\ref{tbl:rexi-ds-gmres-timings}.
We increase the number of cores that we use for the space-parallelization
until the speedup saturates.
This happens at about $192$ cores.
Then we start adding more cores by
increasing the time-parallelization.
Instead of using $384$ cores in space, we compute $8$ summands of the rational approximation in parallel, each using $48$ cores for the matrix-vector operations.
This way, the same problem is solved, but the speedup is increased from $151.15$ to $254.50$.
Doubling the number of cores in time to the maximum of $16$ parallel summands, we get a maximum speedup of $335.15$ on $768$ \ac{CPU} cores.
Note that already when going from $1\times 48$ to $8\times 48$ instead of $1\times 384$ (time $\times$ space) the speedup is better.
We see that by combining the space-parallelization with the time-parallelization we can increase the speedup substantially.

Note that an interesting effect is observed when using $1 \times 96$
cores---the efficiency is larger than one. We assume
that this effect is caused by a better \ac{CPU} utilization. The smaller
the problem per node gets the less data needs to be stored on one node.
Thus, at some point, the whole problem fits into the \ac{CPU} cache of the
node. Adding more (spatial) nodes to the problem, however, degrades efficiency
severely.


\section{Conclusions and Outlook}
\label{sec:outro}

In this paper we have derived and applied a new variant of the rational approximation of exponential integrators (\ac{REXI}) approach for the non-relativistic, single-particle Schrödinger equation.
This time-integration scheme, being already more efficient than the standard integrators used for the examples in this work, can be parallelized efficiently.
Each summand of the approximation can be computed in parallel, thus implementing a parallel-across-the-method approach, which augments a classical parallelization strategy in space.
With this approach, scaling limits of distributed matrix- and vector-operations that correspond to operations in the spatial domain can be overcome.
While parallel-in-time techniques are rather successful for problems of parabolic-type, propagations of waves like in the case of the Schrödinger equation are hard to tackle.
With the REXI variant presented here, solving wave-type problems in a time-parallel manner is indeed possible, making efficient fully space-time parallel simulations of quantum systems with the Schrödinger equation possible for the first time.

We have derived and explained the rational approximation strategy chosen for this problem in detail, making use of the Faber-Cara\-théo\-dory-Fejér approximation to compute the shifts and coefficients of the rational approximation of the matrix exponential.
The derivation of the approximation algorithm in Section~\ref{sec:time} can be used as a single-source reference to reproduce or potentially improve the numerical properties of this integrator.
While the classical \ac{REXI} method \citep{HautEtAl2016} is originally tailored for real-valued problems, this approach is also capable of dealing with complex-valued solutions in an efficient way.
In comparison, fewer summands are necessary to achieve the same accuracy, leading to an improved ratio of accuracy per parallel task.
We have shown along the lines of two challenging, real-world examples the impact of the parallel-in-time integrator, in particular with respect to a standard spatial parallelization technique.

For this work we have exclusively focused on the time-dependent, single-particle Schrödinger equation.
The parallel-in-time method used and extended here was motivated by this equation, but its application is not limited to this particular problem.
The approach can be extended to the many-particle Schrödinger equation and, using Newton's method or a suitable implicit-explicit splitting strategy like spectral deferred corrections \citep{Minion2003}, general nonlinear Schrödinger equations can be addressed.
However, there are features of the spatial discretization scheme, which actually limit the potential speedup gained by the \ac{REXI} approach itself.

When assessing the potential of a parallel method, it is important to compute
the speedup with respect to the fastest serial method available. In the case
we considered in this paper, the \ac{REXI} method was also the fastest serial
method. This, however, is in general not the case.

Let us discuss some different situations in which we compare the \ac{REXI}
and the Chebyshev method, to highlight the factors that need to be taken
into account when determining the speedup the \ac{REXI} method can provide.
We do not consider non-exponential methods like Crank-Nicolson, since there the time-step size is prohibitively small.
For the sake of simplicity, we restrict ourselves to comparing the dominant
costs of both methods.
The dominant cost of the \ac{REXI} method is the
solution of $n_\mathrm{R}$ linear
systems, while the dominant cost of the Chebyshev method is the computation
of $n_\mathrm{C}$ matrix vector products involving the matrix $M$
\eqref{eq:matrix-exp-solution}.
In the case that we considered in Section~\ref{sec:num}, $n_\mathrm{R} = 16$ and $n_\mathrm{C} = 26$.
In general, solving a linear system is much more expensive than computing
a matrix-vector product. The reason why in our case the serial \ac{REXI}
method is faster than the Chebyshev method is that $M = -\iu B^{-1}A$, i.e.,
applying $M$ to a vector involves solving a linear system as well.
If we use a discretization in which $B=I$, e.g., a finite difference
discretization this argument no longer holds.

Consider the case where $B=I$. If the time it takes to solve one linear
system is longer than it takes to compute $n_\mathrm{C}$ matrix-vector
products, the \ac{REXI} method provides no speedup over the Chebshev method
independent of the number of processors that are used.
Note that for spectral methods with suitable domain geometries, the costs for solving a linear system and applying a matrix to a vector are very similar.
Thus, for those discretizations \ac{REXI} can provide speedup, too, and the original papers did indeed focus on those methods.
In the case that we
considered, $n_\mathrm{C} = 26$. When solving the linear system not in spectral space but with a linear solver like \ac{GMRES}, the statement essentially means that each system must be
solved using fewer than $26$ iterations, which is a severe limit on the
number of iterations.

Let us now assume that we are in a situation where $B=I$ and solving one of
the linear systems in \eqref{eq:rational-evaluation} is actually faster
than $n_\mathrm{C}$ matrix-vector products. If we use the \ac{GMRES} method
to solve the linear systems, we can use a method like the shifted \ac{GMRES}
method \citep{FrommerGlaessner1998}, which is able to solve a set of shifted
linear systems at about the same cost as it takes to solve one system.
While this leads to a very efficient method, it leaves no room for any
speedup due to time-parallelization.
If we now assume that we need to precondition the \ac{GMRES} iteration and
each shift required a different preconditioner, we can no longer apply the
shifted \ac{GMRES} method. Hence, in this situation, is is again possible to
obtain a speedup using time-parallelization.

Thus, in the case where solving a linear system involving the
matrix $B$ is expensive enough, using the time-parallelization of the
\ac{REXI} method provides a speedup over solving sequentially.
If solving these linear systems is cheap, it is not clear that the \ac{REXI}
method yields a speedup with respect to a certain sequential method.
In the case of a finite element discretization, we are, however, in the
situation, where solving a linear system is expensive enough to make the
use of the time-parallel \ac{REXI} method beneficial.

\section*{Acknowledgements} 
  The authors would like to thank Martin Gander and Martin Schreiber for their valuable input, in particular during the 8th PinT Workshop at the Center for Interdisciplinary Research in Bielefeld, Germany. The authors furthermore thankfully acknowledge the provision of computing time on the JURECA cluster at Jülich Supercomputing Centre.

\bibliography{references}

\begin{thebibliography}{58}
\expandafter\ifx\csname natexlab\endcsname\relax\def\natexlab#1{#1}\fi
\providecommand{\url}[1]{\texttt{#1}}
\providecommand{\href}[2]{#2}
\providecommand{\path}[1]{#1}
\providecommand{\DOIprefix}{doi:}
\providecommand{\ArXivprefix}{arXiv:}
\providecommand{\URLprefix}{URL: }
\providecommand{\Pubmedprefix}{pmid:}
\providecommand{\doi}[1]{\href{http://dx.doi.org/#1}{\path{#1}}}
\providecommand{\Pubmed}[1]{\href{pmid:#1}{\path{#1}}}
\providecommand{\bibinfo}[2]{#2}
\ifx\xfnm\relax \def\xfnm[#1]{\unskip,\space#1}\fi
\bibitem[{Schrödinger(1926)}]{Schroedinger1926}
\bibinfo{author}{E.~Schrödinger}, \bibinfo{journal}{Ann. Phys.}
  \bibinfo{volume}{384} (\bibinfo{year}{1926}) \bibinfo{pages}{361--376}.
  \DOIprefix\doi{10.1002/andp.19263840404}.
\bibitem[{Mazur and Rubin(1959)}]{MazurRubin1959}
\bibinfo{author}{J.~Mazur}, \bibinfo{author}{R.~J. Rubin}, \bibinfo{journal}{J.
  Chem. Phys.} \bibinfo{volume}{31} (\bibinfo{year}{1959})
  \bibinfo{pages}{1395--1412}. \DOIprefix\doi{10.1063/1.1730605}.
\bibitem[{McCullough and Wyatt(1969)}]{McCulloughWyatt1969}
\bibinfo{author}{E.~A. McCullough}, \bibinfo{author}{R.~E. Wyatt},
  \bibinfo{journal}{J. Chem. Phys.} \bibinfo{volume}{51} (\bibinfo{year}{1969})
  \bibinfo{pages}{1253--1254}. \DOIprefix\doi{10.1063/1.1672133}.
\bibitem[{McCullough and Wyatt(1971)}]{McCulloughWyatt1971a}
\bibinfo{author}{E.~A. McCullough}, \bibinfo{author}{R.~E. Wyatt},
  \bibinfo{journal}{J. Chem. Phys.} \bibinfo{volume}{54} (\bibinfo{year}{1971})
  \bibinfo{pages}{3578--3591}. \DOIprefix\doi{10.1063/1.1675384}.
\bibitem[{Park et~al.(1970)Park, Tahk, and Wilson}]{ParkEtAl1970}
\bibinfo{author}{Y.~R.~L. Park}, \bibinfo{author}{C.~T. Tahk},
  \bibinfo{author}{D.~J. Wilson}, \bibinfo{journal}{J. Chem. Phys.}
  \bibinfo{volume}{53} (\bibinfo{year}{1970}) \bibinfo{pages}{786--791}.
  \DOIprefix\doi{10.1063/1.1674059}.
\bibitem[{Balint-Kurti(2010)}]{BalintKurti2010}
\bibinfo{author}{G.~G. Balint-Kurti}, \bibinfo{journal}{Theor. Chem. Acc.}
  \bibinfo{volume}{127} (\bibinfo{year}{2010}) \bibinfo{pages}{1--17}.
  \DOIprefix\doi{10.1007/s00214-010-0760-4}.
\bibitem[{Li(2019)}]{Li2019}
\bibinfo{author}{X.~Li}, \bibinfo{journal}{J. Chem. Phys.}
  \bibinfo{volume}{150} (\bibinfo{year}{2019}) \bibinfo{pages}{114111}.
  \DOIprefix\doi{10.1063/1.5079326}.
\bibitem[{Ullrich(2016)}]{Ullrich2016}
\bibinfo{author}{C.~A. Ullrich}, \bibinfo{title}{Time-Dependent
  Density-Functional Theory. Concepts and Applications},
  \bibinfo{publisher}{Oxford University Press}, \bibinfo{year}{2016}.
\bibitem[{Messiah(1999)}]{Messiah1999}
\bibinfo{author}{A.~Messiah}, \bibinfo{title}{Quantum Mechanics},
  \bibinfo{publisher}{Dover Publications}, \bibinfo{year}{1999}.
\bibitem[{Shankar(2014)}]{Shankar2014}
\bibinfo{author}{R.~Shankar}, \bibinfo{title}{Principles of Quantum Mechanics},
  \bibinfo{edition}{2nd} ed., \bibinfo{publisher}{Springer},
  \bibinfo{year}{2014}.
\bibitem[{Griffiths(2017)}]{Griffiths2017}
\bibinfo{author}{D.~J. Griffiths}, \bibinfo{title}{Introduction to Quantum
  Mechanics}, \bibinfo{edition}{2nd} ed., \bibinfo{publisher}{Cambridge
  University Press}, \bibinfo{year}{2017}.
  \DOIprefix\doi{10.1017/9781316841136}.
\bibitem[{Burrage(1997)}]{Burrage1997}
\bibinfo{author}{K.~Burrage}, \bibinfo{journal}{Adv. Comput. Math.}
  \bibinfo{volume}{7} (\bibinfo{year}{1997}) \bibinfo{pages}{1--31}.
  \DOIprefix\doi{10.1023/A:1018997130884}.
\bibitem[{Schöbel and Speck(2019)}]{schoebel2019pfasster}
\bibinfo{author}{R.~Schöbel}, \bibinfo{author}{R.~Speck},
  \bibinfo{title}{{PFASST-ER: Combining the Parallel Full Approximation Scheme
  in Space and Time with parallelization across the method}},
  \bibinfo{year}{2019}. \href{http://arxiv.org/abs/1912.00702}{\tt
  arXiv:1912.00702}.
\bibitem[{Clarke et~al.(2020)Clarke, Davies, Ruprecht, and
  Tobias}]{ClarkeEtAl_ParallelInTimeIntegrationOfKinematicDynamos}
\bibinfo{author}{A.~T. Clarke}, \bibinfo{author}{C.~J. Davies},
  \bibinfo{author}{D.~Ruprecht}, \bibinfo{author}{S.~M. Tobias},
  \bibinfo{journal}{J. Comput. Phys. X} \bibinfo{volume}{7}
  (\bibinfo{year}{2020}) \bibinfo{pages}{100057}.
  \DOIprefix\doi{10.1016/j.jcpx.2020.100057}.
\bibitem[{Friedhoff et~al.(2019)Friedhoff, Hahne, and
  Schöps}]{FriedhoffEtAl_MultigridReductionInTimeForEddyCurrentProblems}
\bibinfo{author}{S.~Friedhoff}, \bibinfo{author}{J.~Hahne},
  \bibinfo{author}{S.~Schöps}, \bibinfo{journal}{Proc. Appl. Math. Mech.}
  \bibinfo{volume}{19} (\bibinfo{year}{2019}) \bibinfo{pages}{e201900262}.
  \DOIprefix\doi{10.1002/pamm.201900262}.
\bibitem[{Samaddar et~al.(2019)Samaddar, Coster, Bonnin, Berry, Elwasif, and
  Batchelor}]{SamaddarEtAl2019}
\bibinfo{author}{D.~Samaddar}, \bibinfo{author}{D.~Coster},
  \bibinfo{author}{X.~Bonnin}, \bibinfo{author}{L.~Berry},
  \bibinfo{author}{W.~Elwasif}, \bibinfo{author}{D.~Batchelor},
  \bibinfo{journal}{Comput. Phys. Commun.} \bibinfo{volume}{235}
  (\bibinfo{year}{2019}) \bibinfo{pages}{246--257}.
  \DOIprefix\doi{10.1016/j.cpc.2018.08.007}.
\bibitem[{Schroder et~al.(2018)Schroder, Falgout, Woodward, Top, and
  Lecouvez}]{SchroderEtAl2017}
\bibinfo{author}{J.~B. Schroder}, \bibinfo{author}{R.~D. Falgout},
  \bibinfo{author}{C.~S. Woodward}, \bibinfo{author}{P.~Top},
  \bibinfo{author}{M.~Lecouvez}, in: \bibinfo{booktitle}{2018 IEEE Power \&
  Energy Society General Meeting (PESGM)}, \bibinfo{organization}{IEEE}, pp.
  \bibinfo{pages}{1--5}.
\bibitem[{Agboh et~al.(2019)Agboh, Ruprecht, and Dogar}]{agboh2019combining}
\bibinfo{author}{W.~C. Agboh}, \bibinfo{author}{D.~Ruprecht},
  \bibinfo{author}{M.~R. Dogar}, \bibinfo{title}{Combining coarse and fine
  physics for manipulation using parallel-in-time integration},
  \bibinfo{year}{2019}. \href{http://arxiv.org/abs/1903.08470}{\tt
  arXiv:1903.08470}.
\bibitem[{Trefethen and
  Weideman(2014)}]{TrefethenWeideman_TheExponentiallyConvergentTrapezoidalRule}
\bibinfo{author}{L.~N. Trefethen}, \bibinfo{author}{J.~A.~C. Weideman},
  \bibinfo{journal}{SIAM Rev.} \bibinfo{volume}{56} (\bibinfo{year}{2014})
  \bibinfo{pages}{385--458}. \DOIprefix\doi{10.1137/130932132}.
\bibitem[{Hale et~al.(2008)Hale, Higham, and
  Trefethen}]{HaleEtAl_ComputingPowerLogAndRelatedMatrixFunctionsByContourIntegrals}
\bibinfo{author}{N.~Hale}, \bibinfo{author}{N.~Higham},
  \bibinfo{author}{L.~Trefethen}, \bibinfo{journal}{SIAM J. Numer. Anal.}
  \bibinfo{volume}{46} (\bibinfo{year}{2008}) \bibinfo{pages}{2505--2523}.
  \DOIprefix\doi{10.1137/070700607}.
\bibitem[{Haut et~al.(2016)Haut, Babb, Martinsson, and Wingate}]{HautEtAl2016}
\bibinfo{author}{T.~S. Haut}, \bibinfo{author}{T.~Babb}, \bibinfo{author}{P.~G.
  Martinsson}, \bibinfo{author}{B.~A. Wingate}, \bibinfo{journal}{IMA J. Numer.
  Anal.} \bibinfo{volume}{36} (\bibinfo{year}{2016}) \bibinfo{pages}{688--716}.
  \DOIprefix\doi{10.1093/imanum/drv021}.
\bibitem[{Schreiber et~al.(2019)Schreiber, Schaeffer, and
  Loft}]{SchreiberEtAl_ExponentialIntegratorsWithParallelInTimeRationalApproximationsForTheShallowWaterEquationsOnTheRotatingSphere}
\bibinfo{author}{M.~Schreiber}, \bibinfo{author}{N.~Schaeffer},
  \bibinfo{author}{R.~Loft}, \bibinfo{journal}{Parallel Comput.}
  (\bibinfo{year}{2019}). \DOIprefix\doi{10.1016/j.parco.2019.01.005}.
\bibitem[{Schreiber et~al.(2018)Schreiber, Peixoto, Haut, and
  Wingate}]{SchreiberEtAl2018}
\bibinfo{author}{M.~Schreiber}, \bibinfo{author}{P.~S. Peixoto},
  \bibinfo{author}{T.~Haut}, \bibinfo{author}{B.~Wingate},
  \bibinfo{journal}{Int. J. High Perform. C.} \bibinfo{volume}{32}
  (\bibinfo{year}{2018}) \bibinfo{pages}{913--933}.
  \DOIprefix\doi{10.1177/1094342016687625}.
\bibitem[{Johnson(2009)}]{JohNumerical2009}
\bibinfo{author}{C.~Johnson}, \bibinfo{title}{Numerical Solution of Partial
  Differential Equations by the Finite Element Method},
  \bibinfo{publisher}{Dover Publications}, \bibinfo{year}{2009}.
\bibitem[{Brenner and Scott(2000)}]{BrennerScott2000}
\bibinfo{author}{S.~C. Brenner}, \bibinfo{author}{L.~R. Scott},
  \bibinfo{title}{The Mathematical Theory of Finite Element Methods},
  \bibinfo{publisher}{Springer}, \bibinfo{year}{2000}.
\bibitem[{Braess(2007)}]{Braess2007}
\bibinfo{author}{D.~Braess}, \bibinfo{title}{Finite Elements. Theory, Fast
  Solvers, and Applications in Elasticity Theory}, \bibinfo{edition}{3rd} ed.,
  \bibinfo{publisher}{Cambridge University Press}, \bibinfo{year}{2007}.
\bibitem[{Arai et~al.(1976)Arai, Kanesaka, and Kagawa}]{AraiEtAl1976}
\bibinfo{author}{H.~Arai}, \bibinfo{author}{I.~Kanesaka},
  \bibinfo{author}{Y.~Kagawa}, \bibinfo{journal}{Bull. Chem. Soc. Jpn.}
  \bibinfo{volume}{49} (\bibinfo{year}{1976}) \bibinfo{pages}{1785--1787}.
  \DOIprefix\doi{10.1246/bcsj.49.1785}.
\bibitem[{Kanesaka et~al.(1978)Kanesaka, Arai, and Kawai}]{KanesakaEtAl1978}
\bibinfo{author}{I.~Kanesaka}, \bibinfo{author}{H.~Arai},
  \bibinfo{author}{K.~Kawai}, \bibinfo{journal}{Bull. Chem. Soc. Jpn.}
  \bibinfo{volume}{51} (\bibinfo{year}{1978}) \bibinfo{pages}{28--32}.
  \DOIprefix\doi{10.1246/bcsj.51.28}.
\bibitem[{Ritz(1909)}]{Ritz1909}
\bibinfo{author}{W.~Ritz}, \bibinfo{journal}{J. Reine Angew. Math.}
  \bibinfo{volume}{135} (\bibinfo{year}{1909}) \bibinfo{pages}{1--61}.
  \DOIprefix\doi{10.1515/crll.1909.135.1.}
\bibitem[{Leforestier et~al.(1991)Leforestier, Bisseling, Cerjan, Feit,
  Friesner, Guldberg, Hammerich, Jolicard, Karrlein, Meyer, Lipkin, Roncero,
  and Kosloff}]{LeforestierEtAl1991}
\bibinfo{author}{C.~Leforestier}, \bibinfo{author}{R.~Bisseling},
  \bibinfo{author}{C.~Cerjan}, \bibinfo{author}{M.~Feit},
  \bibinfo{author}{R.~Friesner}, \bibinfo{author}{A.~Guldberg},
  \bibinfo{author}{A.~Hammerich}, \bibinfo{author}{G.~Jolicard},
  \bibinfo{author}{W.~Karrlein}, \bibinfo{author}{H.-D. Meyer},
  \bibinfo{author}{N.~Lipkin}, \bibinfo{author}{O.~Roncero},
  \bibinfo{author}{R.~Kosloff}, \bibinfo{journal}{J. Comput. Phys.}
  \bibinfo{volume}{94} (\bibinfo{year}{1991}) \bibinfo{pages}{59 -- 80}.
  \DOIprefix\doi{10.1016/0021-9991(91)90137-A}.
\bibitem[{Bellman(1997)}]{Bellman1997}
\bibinfo{author}{R.~Bellman}, \bibinfo{title}{Introduction to Matrix Analysis},
  \bibinfo{edition}{2nd} ed., \bibinfo{publisher}{SIAM}, \bibinfo{year}{1997}.
\bibitem[{Liesen and Mehrmann(2015)}]{LiesenMehrmann2015}
\bibinfo{author}{J.~Liesen}, \bibinfo{author}{V.~Mehrmann},
  \bibinfo{title}{Linear Algebra}, Springer Undergraduate Mathematics Series,
  \bibinfo{publisher}{Springer}, \bibinfo{year}{2015}.
  \DOIprefix\doi{10.1007/978-3-319-24346-7}.
\bibitem[{Moler and Van~Loan(1978)}]{MolerVanLoan1978}
\bibinfo{author}{C.~Moler}, \bibinfo{author}{C.~Van~Loan},
  \bibinfo{journal}{SIAM Rev.} \bibinfo{volume}{20} (\bibinfo{year}{1978})
  \bibinfo{pages}{801--836}. \DOIprefix\doi{10.1137/1020098}.
\bibitem[{Moler and Van~Loan(2003)}]{MolerVanLoan2003}
\bibinfo{author}{C.~Moler}, \bibinfo{author}{C.~Van~Loan},
  \bibinfo{journal}{SIAM Rev.} \bibinfo{volume}{45} (\bibinfo{year}{2003})
  \bibinfo{pages}{3--49 (electronic)}. \DOIprefix\doi{10.1137/S00361445024180}.
\bibitem[{Higham(2008)}]{Higham2008}
\bibinfo{author}{N.~J. Higham}, \bibinfo{title}{Functions of Matrices. Theory
  and Computation}, \bibinfo{publisher}{SIAM}, \bibinfo{year}{2008}.
\bibitem[{Ellacott(1983)}]{EllFaber1983}
\bibinfo{author}{S.~Ellacott}, \bibinfo{journal}{SIAM J. Numer. Anal.}
  \bibinfo{volume}{20} (\bibinfo{year}{1983}) \bibinfo{pages}{989--1000}.
  \DOIprefix\doi{10.1137/0720069}.
\bibitem[{Trefethen(1981)}]{Trefethen1981b}
\bibinfo{author}{L.~N. Trefethen}, \bibinfo{journal}{Numer. Math.}
  \bibinfo{volume}{37} (\bibinfo{year}{1981}) \bibinfo{pages}{297--320}.
  \DOIprefix\doi{10.1007/BF01398258}.
\bibitem[{Greene and Kim(2017)}]{GreeneKim2017}
\bibinfo{author}{R.~E. Greene}, \bibinfo{author}{K.-T. Kim},
  \bibinfo{journal}{Complex Anal. Synerg.} \bibinfo{volume}{3}
  (\bibinfo{year}{2017}) \bibinfo{pages}{1}.
  \DOIprefix\doi{10.1186/s40627-016-0009-7}.
\bibitem[{Curtiss(1971)}]{CurFaber1971}
\bibinfo{author}{J.~H. Curtiss}, \bibinfo{journal}{Amer. Math. Monthly}
  \bibinfo{volume}{78} (\bibinfo{year}{1971}) \bibinfo{pages}{577--596}.
\bibitem[{Stewart(2002)}]{Stewart_AKrylovSchurAlgorithmForLargeEigenproblems}
\bibinfo{author}{G.~W. Stewart}, \bibinfo{journal}{SIAM J. Matrix Anal. Appl.}
  \bibinfo{volume}{23} (\bibinfo{year}{2002}) \bibinfo{pages}{601--614}.
  \DOIprefix\doi{10.1137/S0895479800371529}.
\bibitem[{Quarteroni et~al.(2000)Quarteroni, Sacco, and
  Saleri}]{QuarteroniEtAl_NumericalMethematics}
\bibinfo{author}{A.~Quarteroni}, \bibinfo{author}{R.~Sacco},
  \bibinfo{author}{F.~Saleri}, \bibinfo{title}{Numerical Methematics},
  number~\bibinfo{number}{37} in \bibinfo{series}{Text in Applied Mathematics},
  \bibinfo{publisher}{Springer}, \bibinfo{year}{2000}.
\bibitem[{Hairer and
  Wanner(2002)}]{HairerWanner_SolvingOrdinaryDifferentialEquationsII}
\bibinfo{author}{E.~Hairer}, \bibinfo{author}{G.~Wanner},
  \bibinfo{title}{Solving Ordinary Differential Equations II. Stiff and
  Differential-Algebraic Problems}, Springer Series in Computational
  Mathematics, \bibinfo{edition}{2nd} ed., \bibinfo{publisher}{Springer},
  \bibinfo{year}{2002}. \DOIprefix\doi{10.1007/978-3-642-05221-7}.
\bibitem[{Tal‐Ezer and Kosloff(1984)}]{TalEzerKosloff1984}
\bibinfo{author}{H.~Tal‐Ezer}, \bibinfo{author}{R.~Kosloff},
  \bibinfo{journal}{J. Chem. Phys.} \bibinfo{volume}{81} (\bibinfo{year}{1984})
  \bibinfo{pages}{3967--3971}. \DOIprefix\doi{10.1063/1.448136}.
\bibitem[{Kirk et~al.(2006)Kirk, Peterson, Stogner, and Carey}]{KirkEtAl2006}
\bibinfo{author}{B.~S. Kirk}, \bibinfo{author}{J.~W. Peterson},
  \bibinfo{author}{R.~H. Stogner}, \bibinfo{author}{G.~F. Carey},
  \bibinfo{journal}{Eng. Comput.} \bibinfo{volume}{22} (\bibinfo{year}{2006})
  \bibinfo{pages}{237--254}. \DOIprefix\doi{10.1007/s00366-006-0049-3}.
\bibitem[{PETSc(????)}]{Petsc}
PETSc, \bibinfo{title}{{PETSc} website}, ???? \URLprefix
  \url{https://www.mcs.anl.gov/petsc}.
\bibitem[{{Jülich Supercomputing Centre}(2018)}]{JSC2018}
\bibinfo{author}{{Jülich Supercomputing Centre}}, \bibinfo{journal}{J.
  Large-Scale Res. Facilities} \bibinfo{volume}{4} (\bibinfo{year}{2018})
  \bibinfo{pages}{A132}. \DOIprefix\doi{10.17815/jlsrf-4-121-1}.
\bibitem[{Hartree(1928)}]{Hartree1928a}
\bibinfo{author}{D.~R. Hartree}, \bibinfo{journal}{Math. Proc. Cambridge
  Philos. Soc.} \bibinfo{volume}{24} (\bibinfo{year}{1928})
  \bibinfo{pages}{89–110}. \DOIprefix\doi{10.1017/S0305004100011919}.
\bibitem[{Mills et~al.(1993)Mills, Cvitaš, Homann, Kallay, and
  Kuchitsu}]{MillsEtAl1993}
\bibinfo{editor}{I.~Mills}, \bibinfo{editor}{T.~Cvitaš},
  \bibinfo{editor}{K.~Homann}, \bibinfo{editor}{N.~Kallay},
  \bibinfo{editor}{K.~Kuchitsu} (Eds.), \bibinfo{title}{Quantities, Units and
  Symbols in Physical Chemistry}, \bibinfo{edition}{2nd} ed.,
  \bibinfo{publisher}{Blackwell Science}, \bibinfo{year}{1993}.
\bibitem[{Rivlin(1981)}]{Rivlin1981}
\bibinfo{author}{T.~J. Rivlin}, \bibinfo{title}{An Introduction to the
  Approximation of Functions}, \bibinfo{publisher}{Dover Publications},
  \bibinfo{year}{1981}.
\bibitem[{Powell(1981)}]{Powell1981}
\bibinfo{author}{M.~J.~D. Powell}, \bibinfo{title}{Approximation theory and
  methods}, \bibinfo{publisher}{Cambridge University Press},
  \bibinfo{year}{1981}.
\bibitem[{Trefethen(2012)}]{Trefethen2012}
\bibinfo{author}{L.~N. Trefethen}, \bibinfo{title}{Approximation Theory and
  Approximation Practice}, \bibinfo{publisher}{SIAM}, \bibinfo{year}{2012}.
\bibitem[{Clenshaw(1955)}]{Clenshaw1955}
\bibinfo{author}{C.~W. Clenshaw}, \bibinfo{journal}{Math. Comp.}
  \bibinfo{volume}{9} (\bibinfo{year}{1955}) \bibinfo{pages}{pp.118--120}.
  \DOIprefix\doi{10.1090/S0025-5718-1955-0071856-0}.
\bibitem[{Askar and Cakmak(1978)}]{AskarCakmak1978}
\bibinfo{author}{A.~Askar}, \bibinfo{author}{A.~S. Cakmak},
  \bibinfo{journal}{J. Chem. Phys.} \bibinfo{volume}{68} (\bibinfo{year}{1978})
  \bibinfo{pages}{2794--2798}. \DOIprefix\doi{10.1063/1.436072}.
\bibitem[{Li and Demmel(2003)}]{LiDemmel2003}
\bibinfo{author}{X.~S. Li}, \bibinfo{author}{J.~W. Demmel},
  \bibinfo{journal}{ACM Trans. Math. Software} \bibinfo{volume}{29}
  (\bibinfo{year}{2003}) \bibinfo{pages}{110--140}.
  \DOIprefix\doi{10.1145/779359.779361}.
\bibitem[{Tonomura et~al.(1989)Tonomura, Endo, Matsuda, Kawasaki, and
  Ezawa}]{TonomuraEtAl1989}
\bibinfo{author}{A.~Tonomura}, \bibinfo{author}{J.~Endo},
  \bibinfo{author}{T.~Matsuda}, \bibinfo{author}{T.~Kawasaki},
  \bibinfo{author}{H.~Ezawa}, \bibinfo{journal}{Am. J. Phys}
  \bibinfo{volume}{57} (\bibinfo{year}{1989}) \bibinfo{pages}{117--120}.
  \DOIprefix\doi{10.1119/1.16104}.
\bibitem[{Saad and Schultz(1986)}]{SaadSchultz1986}
\bibinfo{author}{Y.~Saad}, \bibinfo{author}{M.~H. Schultz},
  \bibinfo{journal}{SIAM J. Sci. Statist. Comput.} \bibinfo{volume}{7}
  (\bibinfo{year}{1986}) \bibinfo{pages}{856--869}.
  \DOIprefix\doi{10.1137/0907058}.
\bibitem[{Minion(2003)}]{Minion2003}
\bibinfo{author}{M.~L. Minion}, \bibinfo{journal}{Commun. Math. Sci.}
  \bibinfo{volume}{1} (\bibinfo{year}{2003}) \bibinfo{pages}{471--500}.
\bibitem[{Frommer and Glässner(1998)}]{FrommerGlaessner1998}
\bibinfo{author}{A.~Frommer}, \bibinfo{author}{U.~Glässner},
  \bibinfo{journal}{SIAM J. Sci. Comput.} \bibinfo{volume}{19}
  (\bibinfo{year}{1998}) \bibinfo{pages}{15--26}.
  \DOIprefix\doi{10.1137/S1064827596304563}.

\end{thebibliography}

\end{document}